\documentclass[reqno,10pt]{amsart}
\usepackage{amsmath,amsthm,amsfonts,color,graphicx}
\usepackage[latin1]{inputenc}
\usepackage[makeroom]{cancel}
\usepackage{tikz}
\usepackage{pgfplots}
\usepackage{subfigure,multirow,float}
\usetikzlibrary{decorations.pathreplacing}
\usepackage{filecontents}
\pgfplotsset{compat=1.13}
\usepackage{silence}
\definecolor{lightGray}{RGB}{235,235,235}
\definecolor{orange}{RGB}{255,128,0}
\definecolor{ucib}{RGB}{0,36,105}
\definecolor{mygreen}{RGB}{0,128,0}
\definecolor{green}{RGB}{0,128,0}
\definecolor{lightBlue}{RGB}{102,153,204}

\newtheorem{thm}{Theorem}[section]

\newtheorem{rems}[thm]{Remarks}
\newtheorem{rem}[thm]{Remark}

\DeclareMathAlphabet{\mathpzc}{OT1}{pzc}{m}{it}

\begin{document}
\bibliographystyle{plain}

\title[Smooth Selection Embedding Method]{A Novel Optimization
  Approach to Fictitious Domain Methods}

\author{Daniel Agress}
\author{Patrick Guidotti}
\address{University of California, Irvine\\
Department of Mathematics\\
340 Rowland Hall\\
Irvine, CA 92697-3875\\ USA }
\email{dagress@uci.edu  and gpatrick@math.uci.edu}

\begin{abstract}
A new approach to the solution of boundary value problems
within the so-called {\em fictitious domain methods} philosophy is
proposed which avoids well known shortcomings of other fictitious
domain methods, including the need to generate extensions of the
data. The salient feature of the novel method, which we refer to as
SSEM (Smooth Selection Embedding Method), is that it reduces 
the whole boundary value problem to a linear constraint for an
appropriate optimization problem formulated in a larger, simpler set
containing the domain on which the boundary value problem is posed and
which allows for the use of straightforward
discretizations. The
proposed method in essence computes a (discrete) extension of the
solution to the boundary value problem by selecting it as a smooth
element of the complete affine family of solutions of the original
equations now yielding an under-determined problem for an unkown
defined in the whole fictitious domain. The actual regularity of
this extension is determined by that of the analytic solution and the
choice of objective functional. Numerical experiments will demonstrate
that it can be stably used to efficiently deal with non-constant
coefficients, general geometries, and different boundary conditions in
dimensions $d=1,2,3$ and that it produces solutions of tunable (and
high) accuracy. 
\end{abstract}

\subjclass[1991]{}

\maketitle

\section{Introduction}\label{sec:intro}
In this paper an optimization approach is proposed for the resolution
of general boundary value problems within the framework of {\em
  fictitious domain methods} (we include so-called {\em immersed
  boundary methods} in this class). While the ideas and the methods 
readily apply to any boundary value problem, the approach will be
illustrated by means of second order boundary value problems of
type
\begin{equation}\label{bvp}
  \begin{cases} \mathcal{A} u=f &\text{in }\Omega,\\ \mathcal{B} u=g 
  &\text{on }\Gamma=\partial \Omega, \end{cases}
\end{equation}
for an elliptic operator $\mathcal{A}$ such as, e.g., the Laplacian
$-\Delta$, and an admissible boundary operator $\mathcal{B}$ such as,
e.g., the trace $\gamma_\Gamma$ (Dirichlet problem), the unit outer normal
derivative $\partial _\nu$ (Neumann problem), or a combination thereof
(Robin type problem). Such boundary value problems have traditionally
been strongly or weakly (when in divergence form) formulated as
well-posed problems which admit a unique solution (up to a constant
for some boundary conditions). Most numerical methods, reflecting this
approach and viewpoint, are either a direct discretization of the
problem, like in the case of finite difference methods, or the
discretization of a suitable Dirichlet form-based weak formulation of
the problem, like in the case of finite element methods. When the
domain is special, highly accurate spectral discretizations can be
utilized. The former methods come with the heavy burden 
of generating a mesh for the domain (this becomes a  serious
limiting factor when dealing with some problems, like, for instance,
Moving Boundary Problems or in three space dimensions), whereas the
latter are limited by the small number of allowable shapes for
$\Omega$ and lose some of their benefits for non-constant
coefficients operators. Two widely used methods which seek to avoid
these difficulties are known as the {\em fictitious domain method} and
the {\em immersed boundary method}. These techniques, which we refer to
simply as {\em embedding methods}, transplant the problem from the
original domain $\Omega$ to an encompassing simple region, where
straightforward discretizations and solvers can be utilized. The approach
proposed here can be viewed as a novel embedding method, which reduces
the whole boundary value problem to the role of a
linear constraint to an optimization problem for an appropriately
chosen functional defined on the larger domain. The output of
the method will coincide with an approximation of the solution of the
boundary value problem in the domain $\Omega$ and with a smooth
extension of it defined on
$\mathbb{B}$. The degree of 
smoothness will be determined by the data and the chosen
functional. The method has the advantage 
of working for general domains and general data (read, non-constant
coefficients and any type of boundary conditions) while delivering a
paradigm to obtain, in principle, discretizations of any degree of
accuracy. Not least, it allows for straightforward, robust
implementation, by use of either the QR decomposition or the
preconditioned conjugate gradient method (PCG). It 
differs from other embedding methods in that the 
boundary value problem is left unmodified in the extension process to
the larger domain $\mathbb{B}$. In other words, the interior and
boundary equations are simply discretized by means of the new regular
grid in $\Omega$ and on $\partial \Omega$ for a new ``extended''
unknown vector defined on $\mathbb{B}^m$ (a discretization of
$\mathbb{B}$). A solution is then computed by
selecting a smooth element from the affine space of solutions of the
under-determined problem which results from imposing the equations on the
extended vector.

\subsection{Description of the method.}\label{sec:description} 
As the focus of this paper is on a numerical procedure, the method 
will be described at the discrete level. A parallel continuous
formulation as well as an analysis of the method will be addressed
elsewhere. The continuous counterpart, however, does provide insights
that will be exploited later in the paper in the construction of
effective preconditioners for the iterative PCG-based solution of the
derived equations. For this reason some basic properties of the continuous
operators will be mentioned here and there. 

First fix a simple (square or rectangular) domain $\mathbb{B}$ for
which $\overline{\Omega}\subset \mathbb{B}$.  In this paper
$\mathbb{B}$ will chosen to be the periodic box $(-\pi,\pi)^d\subset
\mathbb{R}^d$. Denote by $\mathbb{B}^m$ a regular uniform
discretization of $\overline{\mathbb{B}}$ consisiting of $N_m$
points, where $m$ is the number of discretization
  points along one and each dimension. Replace the continuous 
differential operator by a discrete counterpart $A=A^m$, defined as a
discrete evaluation of $\mathcal{A}$ at grid-points  which lie inside
$\Omega$
$$
 x\in \Omega^m=\Omega\cap\mathbb{B}^m=\{x_k\, |\, k=1,\dots, N_m^{\Omega}\},\:
 N_m^{\Omega}\in \mathbb{N},
 $$
where $A^m$ acts on ``discrete functions'' defined on
$\mathbb{B}^m$. Given a set of points
$$
 \Gamma^m=\{ y_j\, |\, j=1,\dots, N_m^\Gamma\}\subset \Gamma  
$$
it is possible to discretize the boundary condition using any kind of
interpolation and any kind of discrete differentiation (where needed)
based on the grid $\mathbb{B}^m$ and obtain the corresponding discrete
equation $Bu=B^mu^m=g^n$ for the unknown vector $u^m:\mathbb{B}^m\to
\mathbb{R}$ and a discretization $g^m$ of the boundary function $g$,
defined on $\Gamma^m$. In this way the continuous boundary value
problem \eqref{bvp} can be replaced by the discrete under-determined
system given by
\begin{equation}\label{de}
 Cu=C_mu^m=\begin{bmatrix} A\\ B\end{bmatrix} u=
 \begin{bmatrix} A^m \\ B^m\end{bmatrix} u^m = 
 \begin{bmatrix} f^m \\ g^m\end{bmatrix}=b^{m}=b
\end{equation}
where $f^m$ is a discretization of $f$ at grid points in
$\mathbb{B}^m\cap \Omega$. As the notation indicates, we shall often
suppress the superscripts and the indeces to simplify the
notation. Notice that
\begin{equation*}
 u^m\in \mathbb{R}^{N_m},\: f^m\in \mathbb{R}^{N_m^{\Omega}},
 \text{ and }g^m \in \mathbb{R}^{N_m^\Gamma},
\end{equation*}
for $N_m^{\Omega}=\big | \Omega\cap \mathbb{B}^m\big
|=|\Omega^m|$. Clearly it is always ensured that
$N_m^{\Omega}+N^\Gamma _m<N_m$ so that the problem, while 
under-determined, admits solutions. While not strictly necessary, care
is also taken to make sure that all equations in the system are
independent of each other. The reason is numerical conditioning of the
relevant matrices (more later). Now, and in contrast to available
fictitious domain methods, we don't try to extend or modify the
problem to or in the encompassing domain/grid $\mathbb{B}/\mathbb{B}^m$,
but rather try and find ``the best'' among the solutions of the
under-determined problem \eqref{de}. After all, if you use high order
$\mathbb{B}^m$-based discretizations of derivatives and evaluations, the
equations should be sufficient to determine a solution that achieves
their order of accuracy (up to the order allowed by the regularity of the
solution itself, of course). 

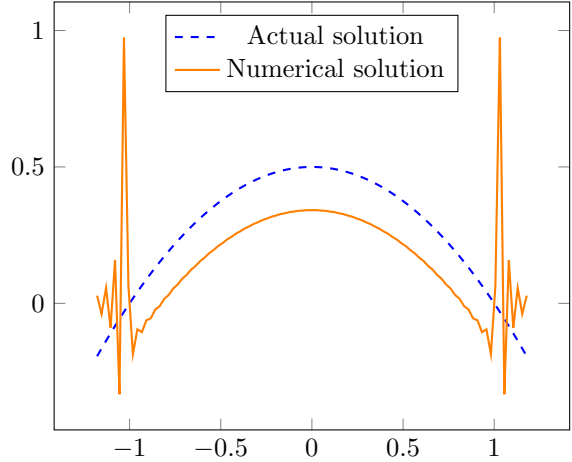
\begin{figure}
  \centering
  \begin{tikzpicture}
	\begin{axis}[legend style={at={(0.5,0.97)},anchor=north}]
	\addplot[mark=none,color=blue,thick,dashed] table{mismatch_sol.data};
	\addlegendentry{Actual solution}
	\addplot[mark=none,color=orange,thick] table{mismatch_u.data};
	\addlegendentry{Numerical solution}
	\end{axis}
  \end{tikzpicture}
  \caption{A 1D visualization of the oscillations caused by trivial
    extension with no regularization. The plot only shows a region that
    is only slightly larger than $\Omega$ since the oscillations 
    occur in a neighborhood of $\partial \Omega$.}
  \label{fig:mismatch}
\end{figure}

A simpleminded approach (which is fine when no regularity at all is
expected) would now be to find a minimal norm solution of the problem,
i.e. solve the linearly constrained optimization problem
\begin{equation}\label{oppb}
  \operatorname{argmin} _{\{Cu=b\}} \frac{1}{2}\| u\| ^2_2,
\end{equation}
where $\|\cdot\| _2$ denotes the Euclidean norm on
$\mathbb{R}^{N_m}$. This would lead to the so-called normal equations 
and to the solution
$$
 u=C^\top \bigl( CC^\top\bigr) ^{-1}b.
$$
Given that the matrix $C=C_{m}$ consists of differential operators
including the evaluation (restriction) in the domain $\Omega^m$ and on
the boundary $\Gamma^m$, its transpose then corresponds to
differential operators containing trivial extensions (read extensions
by 0) and this leads to oscillations generated by the lack of
regularity. This is made apparent in Figure \ref{fig:mismatch}. 
The ``good'' solution is, however, among those of the under-determined
problem, and can be obtained by requiring additional regularity. As
already pointed out, the discretizations $A^m$ and $B^m$ are, after 
all, chosen to have a desired accuracy and the truncations/trivial
extensions destroy it. Thus enforcing an appropriate degree of
regularity should allow for the recovery of the intrinsic accuracy of
the chosen discretizations, again, compatibly with the expected
regularity of the solution itself. This
is also the reason for our choice to call the proposed method Smooth
Selection Embedding Method (SSEM). While this selection is done in a way that is
natural from the point of view of optimization \cite[Chapter
10]{BV04}, it has a nice analytic interpretation which will greatly
help with the practical implementation of the method. Let $\|\cdot\|
_{S}$ be the discretization of a high order norm such as, for instance,
$\| (1-\Delta _\pi)^{p/2}\cdot\|_2$, where $-\Delta_\pi$
denotes the periodic Laplacian on $[-\pi,\pi]^d$ and $p\geq 1$. Now the
problem becomes 
\begin{equation}\label{roppb}
  \operatorname{argmin} _{\{Cu=b\}} \frac{1}{2}\| u\| ^2_{S},
\end{equation}
where the indeces have again been dropped for ease of reading. The
constrained optimization problem \eqref{roppb} can be reformulated as
the unconstrained minimization
$$
 \operatorname{argmin} _{u\in \mathbb{R}^{N_m},\,\Lambda\in
 \mathbb{R}^{N_\Lambda}} \frac{1}{2}\| u\| _{S}^2+\Lambda^\top \bigl(
 Cu-b\bigr), 
$$
upon introduction of Lagrange multipliers $\Lambda\in
\mathbb{R}^{N_\Lambda}$, where $N_\Lambda =N_m^\Omega + N_m^\Gamma$. A
direct computation yields the regularized normal equation
\begin{equation}\label{opsol}
  u=S^{-1}C^\top \bigl( CS^{-1}C^\top\bigr) ^{-1}b,
\end{equation}
where $S$ is the (invertible) operator corresponding to the norm
$\|\cdot\| _{S}$. Now, recalling that $C$ and $C^\top$ are truncated
differential operators (more precisely containing differentiations,
evaluations on subdomains, and extensions), we see that the effect of
the norm is to replace the operator $C^\top$, which, upon being hit by
$C$, is the cause of the oscillations in the simpleminded method, by
the smoothed $S^{-1}C^\top$, which can be captured numerically to a
higher degree of accuracy (no oscillations) when hit by $C$. 
\begin{rem}
While, in the proposed method, Lagrange multipliers are introduced as
they are in many a fictitious domain implementations, the approach is quite
distinct from other methods (see below, Section
\ref{sec:comparison}.) First and foremost the Lagrange 
multipliers are introduced for the whole problem and not
only for the purpose of satisfying the boundary condition. Secondly
they are introduced naturally as an enforcement tool of a linear
constraint and do not require modification of the problem, the use of
extensions, or the introduction of artifical terms (often in the form
of sources).
\end{rem}
\begin{rem}
Notice that formula \eqref{opsol} can be used as a starting point
without any knowledge of a norm generating the operator $S$. One can
choose any convenient smoothing operator acting on (generalized) functions
defined on the box $\mathbb{B}$ instead of $S^{-1}$.
\end{rem}

\subsection{Comparison with Other Embedding Methods} \label{sec:comparison}
Particularly relevant for this paper are the so-called {\em fictitious
  domain methods} and, to a lesser degree {\em immersed boundary
  methods} and {\em boundary integral methods}. These alternative
approaches have experienced a surge in interest in recent years and
seem to be particularly popular in the applied and
very applied communities. Just as with the method advocated here, the
fictitious domain and immersed boundary methods avoid the mesh
generation step by resorting to a ``container'' domain of simple
geometry which admits a straightforward discretization, while boundary
integral methods exploit analytical knowledge about the problem to
obtain a dimensional reduction by
collapsing the problem to the boundary. At the heart of any of these
implementations is the need to resolve the mismatch between the 
boundary and the simple regular grid.
There is a vast literature about these methods as they can be
implemented in various discretization contexts, admit a variety of
distinct practical implementations within each discretization
framework, and can be applied to many
different boundary value problems of mathematical physics \cite{MittalIBM}.
We refer to the beginning of \cite{Low09} for a brief outline of many of these
methods and to \cite{Del03} for a concise description/numerical
implementation of a number of variants. Given the
volume of publications, the choice of references made here was merely
motivated by the fact that they contain a description of the methods'
philosophy and/or many useful additional references in their
introduction. 
\subsubsection{Fictitious Domain Methods} A prominent implementation
procedure, developed by Glowinski and coauthors in
\cite{Glo92,Glo94,Glo07,Glo11} and known as the distributed Lagrange
multiplier method, can be described in some more detail as follows:
think of the domain $\Omega$ as a subset of a larger regular simple
domain $\mathbb{B}$, introduce a (uniform) discretization of
$\mathbb{B}$, and solve the boundary value problem by modifying the
data (the right-hand-side and/or the operator $\mathcal{A}$ in the
prototypical situation considered here), usually by extending them and
by introducing artificially a weighted sum of carefully chosen source
terms supported outside the domain $\Omega$, i.e. in
$\mathbb{B}\setminus \Omega$, or on its boundary $\Gamma$, by
determining the weights (Lagrange multipliers) so as to make sure that
the boundary condition is satisified (or at least
well-approximated). We remark that a common characteristic of these
techniques (and of immersed boundary methods as well) is that Neumann
or Robin boundary conditions are ``natural" and straightforward to
include in the formulation, whereas Dirichlet boundary conditions are
more challenging (see, e.g. \cite{Del03}). These methods clearly have
the advantage of not requiring special care nor effort in the choice
of discretization for $\mathbb{B}$. An often cited criticism of this
approach is the need to extend the original elliptic operator
$\mathcal{A}$ and/or right-hand-side $f$ to corresponding objects
defined on the whole of $\mathbb{B}$. This is not always
straighforward and simple minded extensions (like the trivial one by
zero outside $\Omega$) introduce singularities into the problem
reducing the overall accuracy of the method. See \cite{Boy05}
regarding methods of creating smooth extensions from $\Omega$ to
$\mathbb{B}$ for the purpose of implementing fictitious domain
methods.  
Another approach, in the context of finite elements, consists in
modifying the problem's Dirichlet form to ensure that (non-natural)
boundary conditions be satisified by possibly adding direct or more
subtle penalty or penalty-like terms to it, like, e.g., the so-called
Nitsche method (see \cite{Bur12}, for example). The approach proposed
here can be viewed as a novel fictitious domain method which does not
require any explicit extension of the data (it can itself be used as
remarked later in Section \ref{sec:extension} to produce smooth
extensions) or modification of the original boundary value
problem. Moreover, it makes apparent that the real problem that any
fictitious domain methods has to solve is the selection problem among
the infinitely many solutions of the original problem, which are 
generated as the problem is viewed in a larger domain where it becomes
under-determined. The direct way in which this is done here
(introduction of a high order smoother) clearly shows how the order
of accuracy chosen for the interior and boundary operators can be
recovered in the extended problem through an affine shift obtained by
a natural (both from the point of view of PDEs and of optimization)
regularization.

\subsubsection{Immersed Boundary Methods}
A very popular method used to deal with complex geometries, which is
one of the motivations of this paper as well, is the so-called immersed
boundary method by which a problem is extended to a simple
encompassing domain admitting robust and effective discretizations.
The extension is obtained by the use of Dirac distributions in the distance
from the boundary (more precisely, line and surface integral distributions along
the boundary) and hence typically introduces singularities which
reduce the overall accuracy of the method to first order. Recently, approaches
have been proposed in which the accuracy is improved by the use of extension
operators that preserve smoothness. We refer in particular to
\cite{Ste16} for an immersed boundary method which includes a smooth
extension method, thereby preserving higher order accuracy, albeit at
the cost of significant additional computational time (in what is
called the preparation phase in the paper). We again point out that the
method proposed here does not require any explicit extension since it
identifies the solution among the infinitely many of the extended,
under-determined problem by simply requiring smoothness in the
full computational domain (and hence across the boundary) along with
directly enforcing the PDE in $\Omega$ and the boundary conditions on
$\partial\Omega$ by resorting only to the regular grid.
\subsubsection{Boundary Integral Methods}
While not directly connected to boundary integral methods, the procedure
developed here allows for a nice discrete interpretation of these from
the point of view of optimization. They can be used when 
the existence of an explicit representation for a fundamental solution $G$ of the
differential operator $\mathcal{A}$ is known. In this case one can use
the representation $u_h=\int _\Gamma G(\cdot,y)h(y)\, d \sigma _\Gamma (y)$
for solutions of $\mathcal{A} u=0$ and reduce the boundary value
problem to determining the density $h:\Gamma\to \mathbb{R}$ such that
\begin{equation}\label{bimeq}
  Bu_h(x)=B\int _\Gamma G(x,y) h(y)\, d \sigma _\Gamma (y)=g(x),\:
  x\in \Gamma.
\end{equation}
This effectively leads to a dimensional reduction in the problem as
the unknown density function is only defined on the boundary.

In formulation \eqref{roppb}, this corresponds to situations where the
kernel of $A$ is known and can therefore be represented as the range
of a matrix $M$. In this case, if $u_f$ is a particular solution of
$Au=f$, then the optimization problem can be reduced to 
\begin{equation}\label{bimop}
 \operatorname{argmin} _{\{BMz=g-Bu_f\}}\frac{1}{2}\| u_f+Mz\|_S^2,
\end{equation}
for the unknown (boundary and hence smaller) vector $z$. While the
regularization used here introduces an additional layer not present in a pure
boundary integral formulation, the corresponding problem can also be
efficiently solved given the explicit nature of the smoother and of the 
encompassing domain. Clearly $M$ corresponds to the integral operator
appearing in \eqref{bimeq}, while $B$ is the continuous boundary
operator in \eqref{bimeq} and a corresponding discretization of it in
\eqref{bimop}.
\begin{rems}
We conclude this introduction with a few important remarks.\\
{\bf (a)} The method
is generic in the sense that it does require specific discretizations
of the encompassing domain $\mathbb{B}$ and of the data. It is rather a
procedure that can be adapted to the context of finite differences,
finite elements, or spectral methods quite easily.\\[0.1cm]
{\bf (b)} It has the structure of a classical optimization problem
with linear constraints for which a host of methods exists which can
be used for its resolution. While there seem to be ``natural'' choices
for the smoothing norm $\| \cdot \| _S$, it is possible to work
with other (non-quadratic) functionals, that may deliver better
results for specific problems.\\[0.1cm]
{\bf (c)} It fully avoids the issues related to the need of generating
extensions of the data from the domain $\Omega$ to the encompassing box
$\mathbb{B}$, while, as a matter of fact, it can itself be adapted to produce
smooth extensions. See Section \ref{sec:extension} later in this paper.\\[0.1cm]
{\bf  (d)} As the numerical experiments presented in Section \ref{sec:experiments} will
demonstrate, it is general enough to be robustly implemented for
general domains, for non-constant coefficients, as well as for a
variety of problems (in divergence form and not) and boundary
conditions. In its high order implementations, it clearly heavily
relies on the smoothness of the 
data (and hence of the solution), but can be used for non smooth
problems as well (see Section \ref{sec:nonregular}). Clearly even
better results can be obtained in this case, if specific attention is
paid to the region in which the solution is singular by, e.g.,
introducing a weighted smoothizing norm.
\end{rems} 
\section{Method}\label{sec:method}
\subsection{Methods for Solving the Linear System}\label{sec:system}
Before we turn to describing the actual discretizations of the domain
and the differential operators, we describe two general methods for
solving the linear system in such a way as to obtain a solution of
high accuracy. As described in Section \ref{sec:description}, the boundary value
problem can be 
reduced to finding 
$$
 \operatorname{argmin} _{u\in \mathbb{R}^{N_m},\,\Lambda\in
 \mathbb{R}^{N_\Lambda}} \frac{1}{2}\| u\| _{S_p}^2+\Lambda^\top \bigl(
 Cu-b\bigr).
$$
Here, $||\cdot||^2_{S_p} = ||(1-\Delta_\pi)^{p/2} \cdot||^2$ is a penalty
 norm introduced to enforce the regularity of the
solution across the boundary. Thus, the penalty term imposes the
$H^{p}_\pi(\overline{\mathbb{B}})$ 
regularity of the solution (if at all available; but $p$ can and will
of course be adapted to the solution). The subscript $\pi$ indicates
periodicity. This form of the problem can then be reduced to computing the regularized normal equation
\begin{equation*}
  u=S_p^{-1}C^\top \bigl( CS_p^{-1}C^\top\bigr) ^{-1}b,
\end{equation*}
where the operator $S_p$ is given by
$$S_p u = (1-\Delta_\pi)^p u.$$
A naive approach to solving this linear system would be to directly
invert the (regularized) normal matrix $CS_p^{-1}C^T$. However, such an approach fails to
produce a solution of high accuracy. To obtain such a solution, it is
necessary to either use
\begin{enumerate}
\item[-] a smoother with $p$ very large to obtain a very fast rate
  of convergence, or 
\item[-] a very dense grid, where even a slowly converging solution
  can converge. 
\end{enumerate}
Directly inverting the matrix fails in both of these
approaches. Clearly, for a dense grid, particularly in three
dimensions, it becomes prohibitively expensive to store and directly
invert the normal matrix $CS_p^{-1}C^\top$. On the other hand, the order $p$ of the
smoother can not be pushed too high without hitting the limits of
numerical precision. We recall that the smoother is given by
$(1-\Delta _\pi)^{-p}$. In Fourier space, this corresponds to a
multiplication by the function $(1+|k|^2)^{-p}$. If $k_*$ is largest
mode, as soon $|k|_*^{-2p}$ drops below machine precision, which is
roughly $1e-16$, some matrix entries can no longer be captured
numerically and the benefits of accuracy are lost. For example, on a
grid of size $64^2$, the highest order smoother which can
be used is $p=5$. This greatly limits the accuracy we can obtain.

To remedy these problems, we propose two solutions. The first
continues to use explicit matrices, but uses a QR decomposition to
increase the maximal effective $p$. The second uses an iterative
solver, the PCG method, which relies on an implicit form of the linear
operator rather than an explicit matrix, to allow for solving the
system on larger grids.

\subsubsection{QR Approach}
We consider the QR decomposition of the matrix 
$$S_p^{-1/2}C^\top = QR.$$
Here, using the notation of Section
\ref{sec:description}, $Q \in \mathbb{R}^{ N_m \times N_{\Lambda}}$ is
an orthogonal matrix satisfying $Q^TQ = I$ while $R \in
\mathbb{R}^{N_{\Lambda} \times N_{\Lambda}}$ is an upper triangular
matrix. We then calculate that 
\begin{align*}
S_p^{-1}C^\top(CS_p^{-1}C^\top)^{-1} 
& = S_p^{-1/2}QR(R^\top Q^\top QR)^{-1} \\
& = S_p^{-1/2}QRR^{-1}(R^\top)^{-1} \\
& = S_p^{-1/2}Q(R^\top)^{-1}.
\end{align*}
What makes this method effective is that the
$RR^{-1}$ cancellation reduces the 
power $S_p^{-1}$ in the matrix to $S_p^{-1/2}$. Thus, we are able to
double the order $2p$ of the smoother before the onset of machine
precision limitations. We are therefore able to use any $p \leq 10$
and obtain a highly accurate solution for a coarse grid very
efficiently, as demonstrated by the numerical experiment documented in
Figure \ref{fig:rates1} and Table \ref{tab:times1}. 

\subsubsection{Using the PCG method.}\label{sec:PCG}
An alternative approach consists in using an iterative solver to
deal with the
linear system on very dense grids. As discussed earlier, numerical
limitations will allow us to only use such a method with a smoother of
limited order ($p \leq 4$). However, by 
increasing the grid size, we are able to 
compensate for the smaller order and still obtain an accurate
solution. Because the linear operators $C$, $C^\top$, and $S_p^{-1}$ can be
naturally implemented using sparse matrices (in the case of 
finite difference discretizations) or the FFT (in the case of spectral
discretizations), iterative methods will lend themselves to very fast
computation.  Furthermore, the matrix is positive and symmetric, so a
natural candidate is the 
conjugate gradient method. We note $S_p$ is of order
$-2p$. Thus, it will be very ill-conditioned for large grids and good
preconditioning is necessary. We refer to Section \ref{sec:precond}
for a more detailed description of the preconditioning procedure which
allows for an efficient PCG implementation. We note, however, that the
preconditioning is more effective for the lower order smoothers;
thus, the increased accuracy stemming from the use of
$S_4$ needs to be balanced
against the larger condition 
number, and hence the slower convergence relative to
$S_2$ and $S_3$ for a given grid size. We note that
the $S_2$ has the particular advantage that, with the
preconditioning discussed in Section \ref{sec:precond}, the 
condition number of the corresponding operator remains uniformly
bounded, regardless of grid size. This is because the operator $C$ is
second order in the interior. As $S_2^{-1}$ is of order $-4$,
$CS_2^{-1}C^T$ is of order $0$ in the interior. In Table
\ref{tab:cond1}, we show the growth in condition number for $S_2,
S_3,$ and $S_4$ for the discretization of the disc problem described
in Section \ref{sec:domain}.

\subsubsection{Rates of convergence and contrasting the methods.} \label{sec:contrast}
\begin{figure}[H]
\includegraphics[scale=.45]{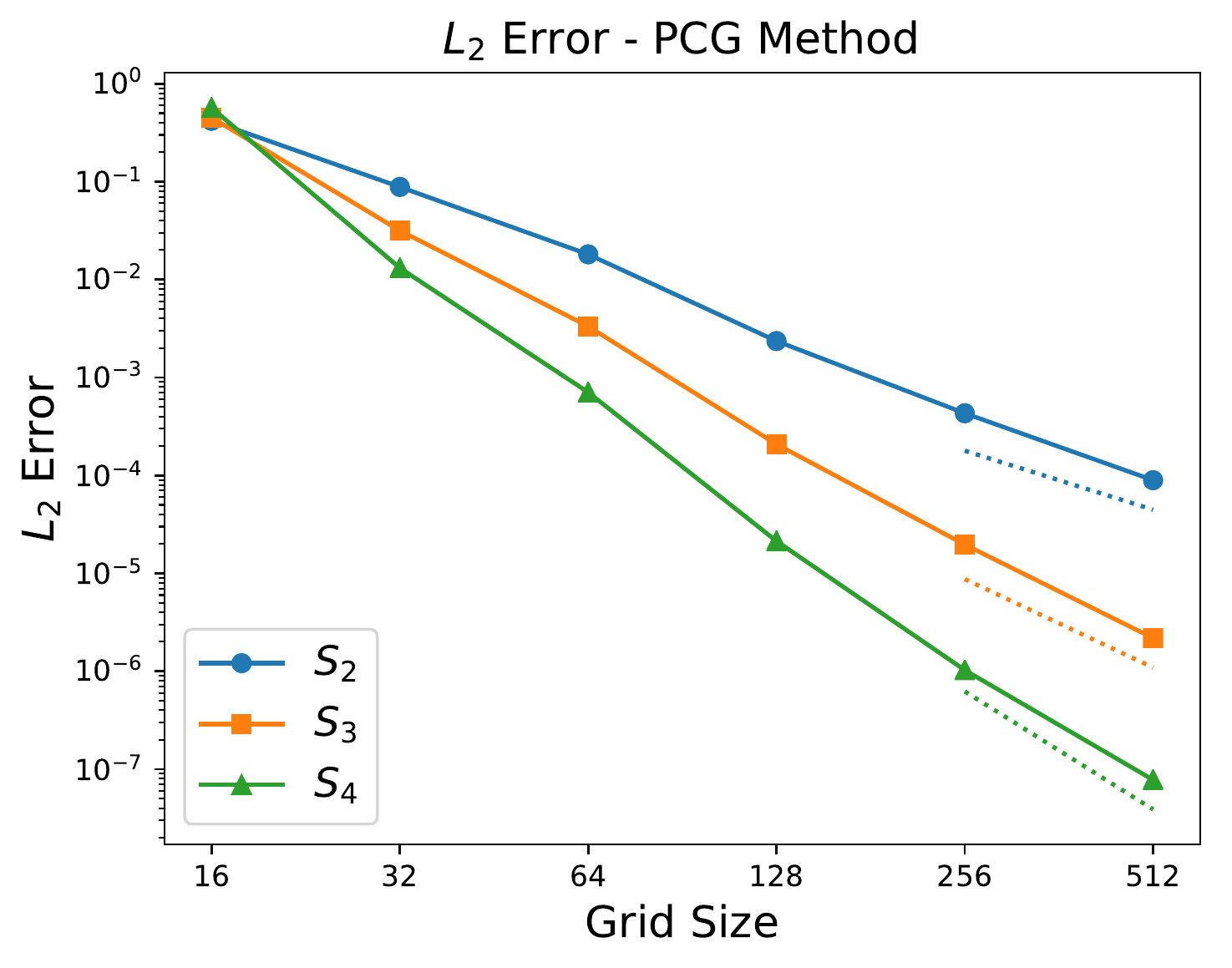}
\includegraphics[scale=.45]{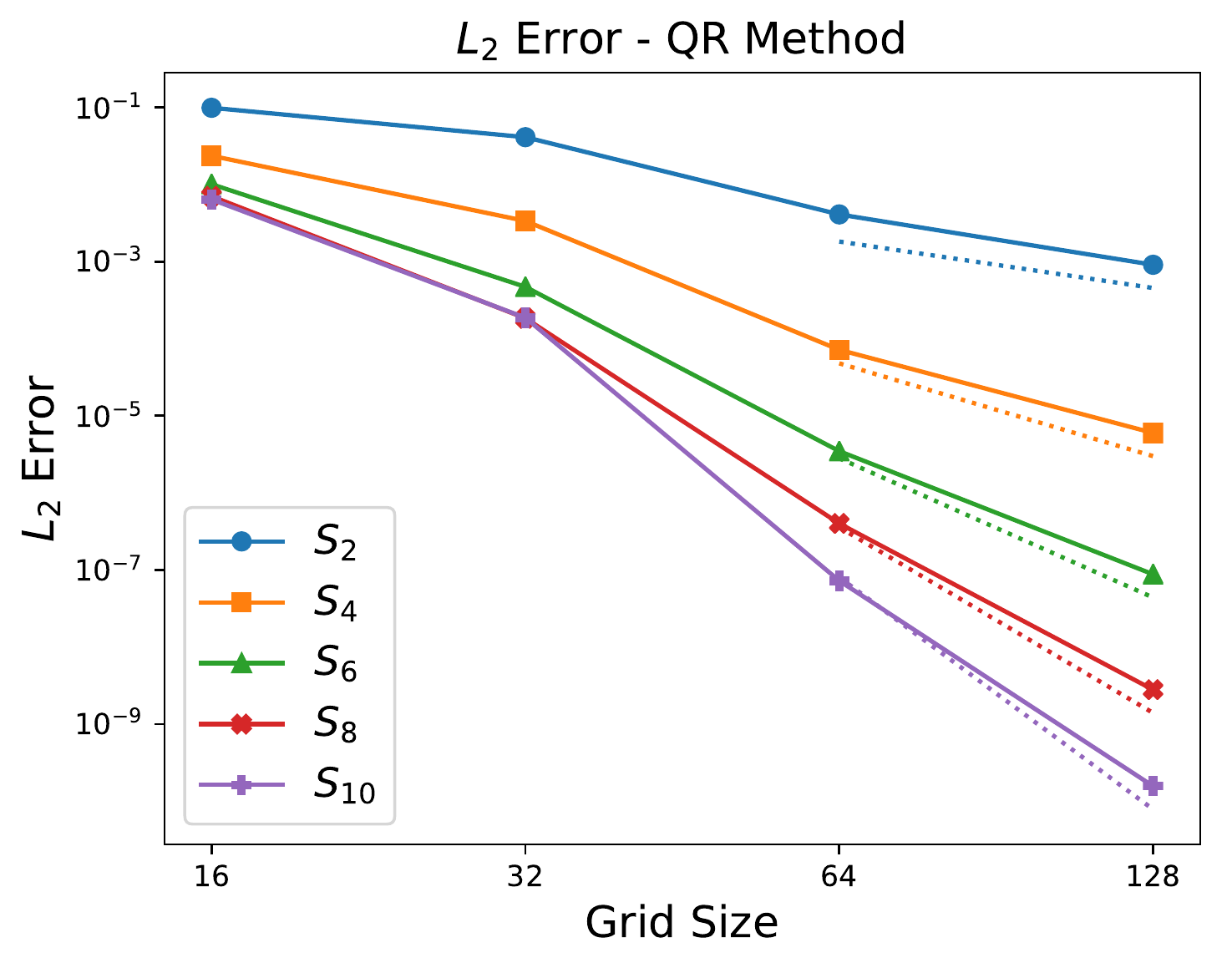}
\caption{Convergence of the $L_2$ error for different order smoothers solving Equation \protect \ref{eq:exp1}. The light dotted lines are reference lines of slope $\frac{1}{m^p}$ where $m$ is the number of grid points along one dimension.}
\label{fig:rates1}
\end{figure}

Next we describe the effectiveness of each of the QR and PCG methods, and
discuss when each should be used. As described in the introduction, a
smoother $S_p$ seeks to find an $H^{p}_\pi(\overline{\mathbb{B}})$ extension of
the solution. Thus, whenever 
the true solution is smooth, we expect that the rate of
convergence of the discrete solution will be of order $p$. In Figure
\ref{fig:rates1}, we demonstrate the rate of convergence of various order
smoothers using both the QR and PCG methods. The 
problem studied is posed on the disc $D$ of radius 1 and reads 
\begin{equation}\label{eq:exp1}
\begin{cases}
-\Delta u = -3(x-y) &\text{ in }D,\\
\phantom{-\Delta}u = x^3 - y^3&\text{ on }\partial D.
\end{cases}
\end{equation}
The exact solution is $x^3 - y^3$. The figure clearly shows the $p$ rate of convergence for
each smoother. We note that the $L_{\infty}$ error converges similarly. Clearly, for a smooth problem, the higher order QR
method on a course grid outerperforms the PCG method, even on a denser
grid. However, in less favorable cases, the PCG method may be
advantageous. For example, for very irregular boundaries, a dense grid
may be necessary to resolve their geometry and the PCG method may be
necessary. Similarly, if the solution itself is not regular, the
higher order smoothers will not achieve faster convergence and it
may be necessary to use the PCG method on a denser grid. Notice,
however, that use of an SVD decomposition on a dense grid is still
possible by using a library which accepts an implicit
linear operator rather than an explicit matrix as its input. This would be an alternative which preserves the accuracy of the QR method with the larger grid of implicit methods. 

\begin{table}[H]
\centering
\begin{tabular}{|c|c|c|c|lllll}
\hline
\multicolumn{1}{|c|}{\multirow{2}{*}{\begin{tabular}[c]{@{}c@{}}Grid  \\ Size\end{tabular}}} & \multicolumn{3}{c|}{CPU Times - PCG Method} & \multicolumn{5}{c|}{CPU Times - QR Method} \\ \cline{2-9} 
\multicolumn{1}{|c|}{} & \multicolumn{1}{c|}{$S_2$} & \multicolumn{1}{c|}{$S_3$} & \multicolumn{1}{c|}{$S_4$} & \multicolumn{1}{c|}{$S_2$} & \multicolumn{1}{c|}{$S_4$} & \multicolumn{1}{c|}{$S_6$} & \multicolumn{1}{c|}{$S_8$} & \multicolumn{1}{c|}{$S_{10}$} \\ \hline
$16^2 $ & $0.01$ & $0.01$ & $0.01$ & \multicolumn{1}{c|}{$0.15$} & \multicolumn{1}{c|}{$0.19$} & \multicolumn{1}{c|}{$0.17$} & \multicolumn{1}{c|}{$0.21$} & \multicolumn{1}{c|}{$0.18$} \\
\hline
$32^2 $ & $0.01$ & $0.02$ & $0.03$ & \multicolumn{1}{c|}{$0.17$} & \multicolumn{1}{c|}{$0.23$} & \multicolumn{1}{c|}{$0.2$} & \multicolumn{1}{c|}{$0.24$} & \multicolumn{1}{c|}{$0.23$} \\
\hline
$64^2 $ & $0.01$ & $0.03$ & $0.07$ & \multicolumn{1}{c|}{$0.39$} & \multicolumn{1}{c|}{$0.56$} & \multicolumn{1}{c|}{$0.52$} & \multicolumn{1}{c|}{$0.56$} & \multicolumn{1}{c|}{$0.61$} \\
\hline
$128^2 $ & $0.04$ & $0.13$ & $0.54$ & \multicolumn{1}{c|}{$6.61$} & \multicolumn{1}{c|}{$7.96$} & \multicolumn{1}{c|}{$8.01$} & \multicolumn{1}{c|}{$8.62$} & \multicolumn{1}{c|}{$7.95$} \\
\hline
$256^2 $ & $0.33$ & $0.59$ & $2.75$& \multicolumn{5}{l}{\multirow{2}{*}{}} \\ \cline{1-4} 
$512^2 $ & $2.15$ & $4.32$ & $24.07$& \multicolumn{5}{l}{} \\ \cline{1-4}
\end{tabular}
\caption{CPU times for solving Equation \protect\ref{eq:exp1}. All
  computations were performed on an Intel 7700HQ.}
\label{tab:times1}
\end{table}

\subsection{Discretization of the Domain}\label{sec:domain}
As described in the introduction, we begin by embedding the domain 
$\Omega$ into a torus $\mathbb{B}$ in order to make use of spectral
methods and of the Fourier transform. The periodicity box $\mathbb{B}$
is discretized with a uniform grid $\mathbb{B}^m$. The boundary $\Gamma$ 
is approximated with a discretization $\Gamma^m$, which is just a set
of $N_m^\Gamma$ points lying on $\partial \Omega$. In practice, it is best for
these points to be uniformly distributed across the boundary. In two
dimensions, this can be accomplished easily by equally spacing points
along an arc length parametrization of the curve. In three dimensions,
equally distributing the points around a surface is more challenging,
although well known algorithms exists for placing points on $S^2$. In
Section \ref{sec:sphere}, we use the well known Fibonacci algorithm
(see \cite{FibSphere}) to create a discretization.

\begin{figure}
\includegraphics[scale=.5]{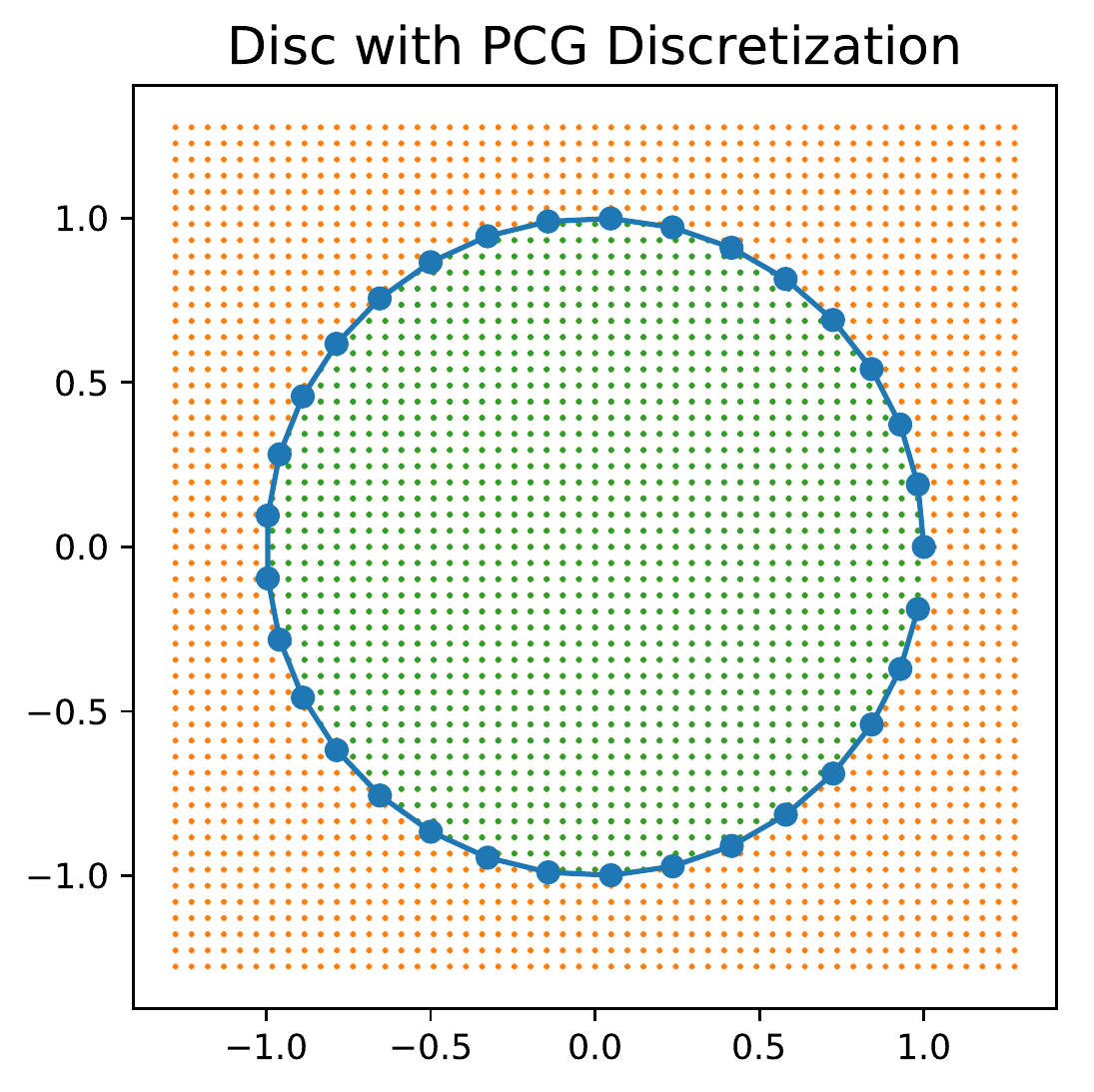}
\includegraphics[scale=.5]{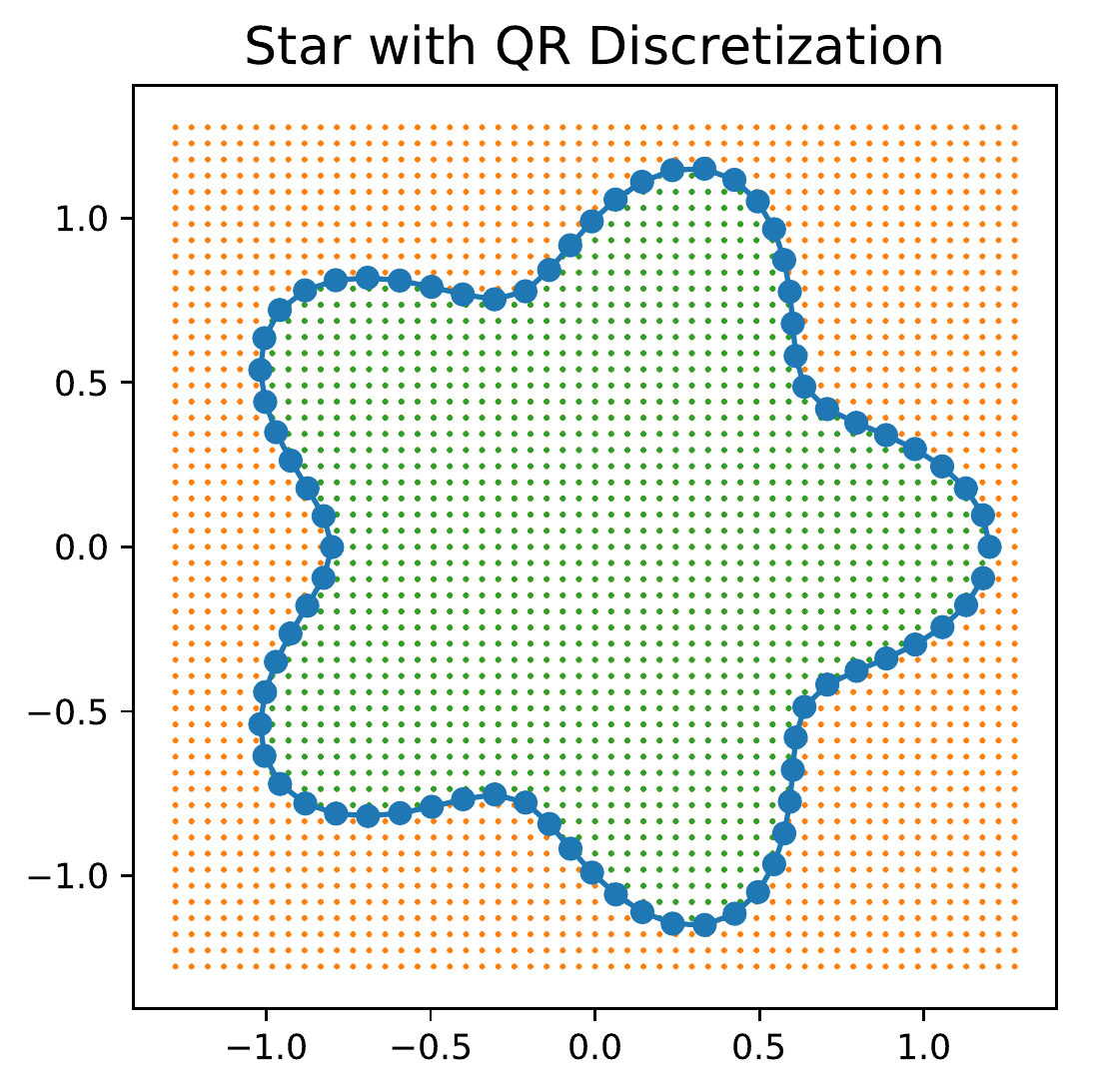}
\caption{Discretizing the boundary.}
\label{fig:domains}
\end{figure}

A choice also needs to be made concerning the density of boundary
points, that is, the value of $n$. When using an insufficient number
of points on the boundary, the accuracy suffers, while too many
points can drive up the condition number. When using
the QR implementation, the method is relatively immune to ill
conditioning, since explicit matrices are used. Thus the boundary
points can be placed close together. If $m$ is the number of grid
points along one dimension, a density of $\frac{1}{2} \frac{m}{2\pi}$
boundary points per unit length seems to be 
effective. The PCG iterative method, on the other hand, is quite
sensitive to ill-conditioning of the matrix. It turns out to be more
effective to space the points further apart according to a density of
$\frac{1}{4}\frac{m}{2\pi}$ points per unit length. This guarantees
that three to four regular grid points lie between any two
boundary points and thereby allows the regular grid $\mathbb{B}^m$
to easily "distinguish" the different boundary points, thereby
keeping the condition number relatively low. In Figure
\ref{fig:domains}, we show the discretization of a disc $D$ with the 
first density described, and a star shaped domain with the second.
For better visualization, we have only plotted the region $[-1.3,1.3]^2$, as
opposed to the entire region $[-\pi,\pi]^2$. In three dimensional
problems, we have found that with a grid of size $m^3$ points, a
boundary spacing of $2\left(\frac{m}{2\pi}\right)^2 $ per unit area
for the QR method is most effective, while
$\frac{1}{16}\left(\frac{m}{2\pi}\right)^2 $ per unit area is best for
the PCG method. This smaller density maintains three 
to four box discretization points between each boundary point along each
dimension, allowing the regular box grid to resolve the ``irregular''
boundary discretization grid. 

\subsection{Discretizing the Differential Operators}
We now discuss the discretization of the differential operators
$C$,$C^\top$, and $S_p$. 
\subsubsection{Discretization of $C$ and $C^\top$.} We recall that 
$$C = \left( \begin{array}{c}
A^m \\ 
B^m
\end{array} \right),$$
where $A^m$ is a matrix of evaluations of a second order differential
operator at the points found in the set $\Omega^m = \mathbb{B}^m\cap
\Omega$, and $B^m$ is a matrix of evaluations of a boundary operator
on the finite subset $\Gamma^m$ of $\partial \Omega$. We begin by
evaluating any necessary derivatives on the entire domain
$\mathbb{B}^m$. As discussed in the introduction, the purpose of using
a fictitious domain method is that it allows us to easily use
techniques which apply to the torus, and extend them to problems with
more complex geometries. In particular, the partial derivatives can be
calculated using either finite difference methods or spectral methods
on the torus. Spectral methods have the advantage of delivering greater
accuracy for smooth problems, while finite difference methods have the
advantage of being slightly faster and being more readily applicable to a wider
range of differential operators. Once the partial derivatives have
been calculated, we restrict the results to $\Omega^m$ and multiply by
the coefficients of the operator $A^m$.

In all of the numerical experiments below, we evaluate the derivatives
used for the operator $A^m$ spectrally. More specifically, whenever
taking the Laplacian, we compute
$$
 (-\Delta)^m
 =\bigl(\mathcal{F}^m\bigr)^{-1}\operatorname{diag}\Bigl(\bigl( 
 |k|^2\bigr) _{k\in \mathbb{Z}^d_m}\Bigr)\mathcal{F}^m,
$$
where $\mathcal{F}^m$ is the discrete fast Fourier transform and $k
\in \mathbb{Z}^d_m$ is the frequency vector at discretization level
$m$. In Subsection \ref{nonconstant}, where we examine 
nonconstant coefficients, we similarly use the Fourier
transform to evaluate the second derivatives in each
  combination of directions. However, we would like to reiterate that $A^m$ 
can implemented with any numerical scheme for calculating derivatives
on the torus. The choices we made were simply dictated by convenience.
Spectral methods are used because we wish to demonstrate the high
order of accuracy which can be obtained by the proposed method.

When applying $(A^m)^\top$, we begin by multiplying by the coefficient of
$A^m$ and, then, take the transpose of the restriction operator part
of $A^m$, which amounts to an extension by $0$ outside of $\Omega^m$.
In this way, we are able to use the chosen method to evaluate the derivatives. 

Because the boundary points $\Gamma^m$ do not lie on the regular grid,
we need to use interpolation operators when implementing the boundary
operator $B^m$. Given that we are interpolating from a regular
rectangular grid, the interpolation operators are simple. Linear,
cubic, or spectral interpolation can all be used. In the examples
below, we have used spectral interpolation. This is because, as
mentioned earlier, we wish to demonstrate the high order of
convergence of the method.

\begin{rem}
The rate of convergence of the solution is constrained by the order of
the smoother, the interpolation operators, and the differential
operators. To avoid wasting computational resources, the order of
accuracy of these various discretizations should be made to match. If
the expected regularity of the solution is known, it can also be taken
into consideration when making this choice.
\end{rem}

\subsubsection{Discretizing the Smoother $S_p$}\label{sec:regularizer}
We now discuss the discretization of the smoother 
$$S_pu = (1-\Delta_{\pi})^p u.$$ 
Because the operator $S_p$ is defined over the torus
$\overline{\mathbb{B}}$, we are able to use the fast Fourier transform
to calculate $S_p^{-1}$, or, as when using the QR method, $S_p^{-1/2}$.
We define the matrix $\mathcal{S}_p$ with diagonal entries
$$\bigl(\mathcal{S}_p\bigr)_{kk}= (1+|k|^2)^{p},$$
where $k \in \mathbb{Z}^d_m$ is the vector of frequencies. We then
note that
$$S_p^{-1}b = \bigl(\mathcal{F}^m\bigr)^{-1} \mathcal{S}_p^{-1} \mathcal{F}^mb.$$
Using the fast Fourier transform, this operator can be evaluated
efficiently with minimal memory requirements. 

\subsubsection{Calculating the Explicit Matrices}
When using the PCG method, the matrix multiplication can be evaluated
implicitly and there is no need to explicitly calculate the matrix
entries. However, the QR decomposition requires an explicit
matrix representation for $C^\top S_p^{-1/2}$. In our implementation, we
have used the simplest option of 
generating the matrix columns by column by evaluating $C^\top S_p^{-1/2}e_i$
for $1\leq i \leq N_{\Lambda}$ for the natural basis vectors
$e_i$. Although this entails many evaluations 
of the matrix, for sparse grids, this time cost is small and the
method is still very efficient. 

An alternative method would be to exploit the fact that the derivative
operators and the smoothing operators are cyclic on the regular
grid. Thus, the matrix can easily be calculated by simply shifting,
for example, $S_p^{-1/2}e_1$ around the grid. A drawback, however, this
method entails explicitly calculating the larger matrix $S_p^{-1/2}$ which
can use large amounts of RAM. We emphasize that the times quoted in
the tables include the time required to calculate the explicit
matrices. We also point out again that libraries exist which can take
the SVD decomposition implicitly; while using the SVD decomposition is
slower than the QR decomposition, doing so would eliminate the need
for the evaluation step.

\subsection{Preconditioning and PCG Implementation}\label{sec:precond}
We now return to a more detailed discussion of the implementation of
the PCG method. As discussed in \ref{sec:PCG}, the normal matrix $CS_p^{-1}C^\top$
is very ill-conditioned and requires a good preconditioner to be
inverted using iterative methods. We note that the ill-conditioning
occurs because of the high order of the operator and because the
boundary operator and the interior operator have different 
orders. To demonstrate this, we think of the operator $CS_p^{-1}C^\top$ as
a block matrix
$$CS_p^{-1}C^\top = \left(\begin{array}{cc}
A^mS_p^{-1}(A^m)^\top & A^m S_p^{-1}(B^m)^\top \\ 
B^n S_p^{-1}(A^m)^\top & B^m S_p^{-1} (B^m)^\top
\end{array} 
\right) = 
\left(\begin{array}{cc}
C_1 & C_2 \\ 
C_2^T & C_3
\end{array} 
\right).$$
As $S_p^{-1}$ is an operator of order $-2p$, the matrix $C_1$ is of
order $4-2p$, $C_2$ is of order $2-2p$ and $C_3$ is of order $-2p$ (for a
boundary operator of order 0). In general, if an operator is of order
$-2p$, the condition number of its matrix will grow like a 
polynomial of degree $2p$ as the 
grid size increases (for example, on
a grid of size $m$, the largest 
eigenvalue of the Laplace operator will be of size $m^2$). Thus, the
large order together with the mismatch in scaling causes a very large 
condition number. We will describe a simple preconditioner which works
effectively for $S_2$, $S_3$, and $S_4$. The
preconditioning consists of finding approximate inverses to the $C_1$
and $C_3$ blocks independently. The general philosophy consists in
preconditioning the operator so that it becomes order $0$. 

We begin by finding an approximate inverse for the $C_3$
block. In the following description, we will consider a Dirichlet
problem, where the boundary operator $\mathcal{B}$ consists of
evaluation on the boundary. The discrete boundary points belonging to
$\Gamma^m$ will be denoted by $y_i$ for $1 \leq i \leq N_m^\Gamma$. We 
recall that
$$B^m: \mathbb{R}^{\mathbb{B}^m} \rightarrow \mathbb{R}^{\Gamma^m}
\text{ and } S_p:\mathbb{R}^{\mathbb{B}^m} \rightarrow
\mathbb{R}^{\mathbb{B}^m}.$$ 
We now consider the operators
$$\widetilde{B}^m: C(\mathbb{B}) \rightarrow \mathbb{R}^{\Gamma^m}
\text{ where } [\widetilde{B}^m u]_i = \langle \delta_{y_i}, u\rangle=u(y_i)$$ 
and 
$$\widetilde{S}_p:H^{2p-d/2-\varepsilon}_\pi(\overline{\mathbb{B}}) \rightarrow
\operatorname{H}^{-d/2-\varepsilon}_\pi(\mathbb{B}) \text{ where }\widetilde{S}_pu = 
(1-\Delta_\pi)^{p}u.$$
We note that $\widetilde{B}^m$ and $\widetilde{S}_p$ can be viewed as
approximations of $B^m$ and $S_p$ respectively, operating on the
continuous $\mathbb{B}$ rather than the discrete $\mathbb{B}^m$. The
integral operator
$$\widetilde{C}_3 := \widetilde{B}^m \widetilde{S}_p^{-1}
(\widetilde{B}^m)^T:\mathbb{R}^{\Gamma^m} \rightarrow
\mathbb{R}^{\Gamma^m} \text{ with kernel } [\widetilde{C}_3]_{ij} =
\langle \delta_{y_i}, \widetilde{S}_p^{-1} \delta_{y_j} \rangle,\: 1\leq
i,j\leq N_m^\Gamma,$$ 
is then a good approximation of $C_3$. Notice that $\delta _y\in
\operatorname{H}^{-d/2-\varepsilon}_\pi(\mathbb{B})$ for any $y\in
\mathbb{B}$ and $\varepsilon>0$. If we define 
$$h(y) = \left( \widetilde{S}_p^{-1} \delta \right)(y)$$
as the fundamental solution of $\widetilde{S}_p$ on the torus
$\mathbb{B}$, we find, by translation invariance of the torus, that
$$[\widetilde{C}_3]_{ij} = \left( \frac{2\pi}{m}\right)^d h(y_i - y_j).$$
Here, by an abuse of notation, the factor $\left(
  \frac{2\pi}{m}\right)^d$ is built-in to 
account for the fact that the ``matrix'' $\widetilde{C}_3$ acts as
an integral operator and not as simply matrix-vector multiplication.
Given a good value table for $h$, we can easily calculate
the matrix $\widetilde{C}_3$ by evaluating the function $h$ on
the matrix of differences between the points in $\Gamma^m$. Given the
explicit matrix $\widetilde{C}_3$, we can directly calculate
$(\widetilde{C}_3)^{-1}$ and use it as a preconditioner for the
$C_3$ block of the matrix. Although this entails inverting a dense
matrix, for coarse grids in three dimensions and even for very fine
grids in two dimensions, the number of boundary points is small enough
that inverting, storing, and applying the matrix is computationally
negligible.

To calculate the function $h$, several methods can be used. In our
implementation, we proceed as follows. We take $\widetilde m$ large
and generate a very fine grid of size $\widetilde m^d$ on the torus
$\overline{\mathbb{B}}$. In our examples, we used $\widetilde m =
4096$. We define the vector $\delta^{\widetilde m}$ by 
$$\delta^{\widetilde m}_k = \begin{cases}
\frac{\widetilde m^d}{(2\pi)^d}, & \text{if }k=0, \\
0, & \text{otherwise.} 
\end{cases}$$
The vector $\delta^{\widetilde m}$ is then an approximation of the
continuous (periodic) $\delta$ distribution supported in the
origin. We then compute $S_p^{-1}\delta^{\widetilde{m}}$ on the fine 
grid. This function is a good approximation of $h$ evaluated at the
points $\mathbb{B}^{\widetilde m}$. We use cubic interpolation to evaluate $h$ at
points which do not lie in $\mathbb{B}^{\widetilde m}$. 
In order to reduce RAM requirements, we only store the numerical
values of $h$ computed by means of the $4096^2$ ($512^3$ in
dimension 3) on a smaller
$256^d$ grid. It is also beneficial to store these values in
memory so they do not need to be recalculated for each problem.  

For the Neumann problem, we note that the order of the matrix $C_3$ is
decreased by $2$, because $C$ and $C^\top$ both evaluate one derivative
on the boundary. Thus, rather than using the function $S_p^{-1}\delta$,
we instead use the function $S^{-1}_{p-1}\delta$. See
Section \ref{sec:neumann} for the effect of this preconditioning for
the Neumann problem.

We now turn to finding an approximate inverse to $C_1$. The matrix
$C_1$ depends on the order of the smoother we have chosen. For $S_2$, we note that the matrix $C_1$ is of order
$0$. Thus, no preconditioning is necessary, and $\widetilde{C}_1^{-1}$
can be taken as the identity. For the$S_3$ and $S_4$, we
note that the operator $C_1$ is the discretization of a differential
operator of order $4-2p$. We wish 
to precondition in such a way as to reduce the order of the operator
to order $0$. Thus, we define the preconditioner
$$\widetilde{C}_1^{-1}u = \bigl(1-\Delta_\Omega\bigr)^{\frac{2p-4}{2}}u.$$
Here, $\Delta_\Omega$ is the Laplace operator on
$\Omega$. In order to implement it, we use the domain discretization
$\Omega^m = \mathbb{B}^m \cap \Omega$ and a finite difference scheme
to discretize the Laplacian on $\Omega^m$. In the examples the five
points stencil (seven points in three dimensions) was chosen to take the
Laplacian on $\Omega^m$. With this preconditioner
$$\widetilde{C}^{-1} = \left( \begin{array}{cc}
\widetilde{C}_1^{-1} & 0 \\ 
0 & \widetilde{C}_3^{-1}
\end{array} \right),$$
the condition number of the preconditioned normal matrix
$\widetilde{C}^{-1/2}(CS_p^{-1}C^\top)\widetilde{C}^{-1/2}$ stays
uniformly bounded, independent of grid size when $p=2$. When $p=3$,
its condition number grows slightly with grid size, while when $p=4$,
it grows significantly with grid 
size; we refer to Table \ref{tab:cond1}. Despite this growth, however, the method is
quite efficient; see Table \ref{tab:times1} as well as the experiments
in Section \ref{sec:experiments} for CPU times. 

\begin{table}[H]\label{tab:cond1}
\centering
\begin{tabular}{|c|c|c|c|c|c|c|c|c|}
\hline 
\multirow{2}{*}{Grid Points} & \multirow{2}{*}{Boundary Points} &
\multicolumn{3}{|c|}{Condition Number} & \multicolumn{3}{|c|}{PCG Iterations} \\  
\cline{3-8}
& & $S_2$ & $S_3$ & $S_4$ & $S_2$ & $S_3$ & $S_4$\\
\hline
$16 \times 16$ & $5$ & $4$ & $17$ & $61$ & $12$ & $20$ & $24$ \\
\hline
$32 \times 32$ & $9$ & $6$ & $20$ & $80$ & $17$ & $30$ & $46$ \\
\hline
$64 \times 64$ & $17$ & $6$ & $23$ & $151$ & $19$ & $33$ & $62$ \\
\hline
$128 \times 128$ & $33$ & $8$ & $25$ & $303$ & $20$ & $37$ & $81$ \\
\hline
$256 \times 256$ & $65$ & $7$ & $25$ & $596$ & $20$ & $33$ & $121$ \\
\hline
\end{tabular}
\caption{Condition numbers and number of iterations for the PCG method
  solving Equation \protect \ref{eq:exp1}. The condition number is
  that of the preconditioned normal matrix $\widetilde{C}^{-1/2}(CS_p^{-1}C^\top)\widetilde{C}^{-1/2}$.} 
\end{table}

\section{Extension Problems}\label{sec:extension}

As described in the introduction, a common problem that is encountered
when embedding a problem with complex geometry in a container space 
is the (smooth) extension of some or all of the data from the original
domain to the encompassing one. We briefly outline how the proposed 
method can be used to generate periodic
$H^{p}_\pi(\overline{\mathbb{B}})$ extensions to the
torus. Other boundary conditions on $\mathbb{B}$ can also be used, of
course, with the appropriate modifications. In this paper we stick to
the periodic setting.

The extension problem consists of finding $\widetilde{u} \in
H^{p}_\pi(\overline{\mathbb{B}})$, given a domain
$\Omega\subset \mathbb{B}$ as  well as a function $u \in
H^{p}(\overline{\Omega})$, satisfying $\widetilde{u} \big|_{\Omega} =
u$. Following the spirit of the 
proposed method, we begin by taking a regular discretization
$\mathbb{B}^m$ of $\mathbb{B}$. We then define 
the operator $C$ to be the restriction to $\Omega^m = \Omega \cap
\mathbb{B}^m$. Correspondingly $C^\top$ is simply given by the extension 
by $0$ from $\Omega^m$ to $\mathbb{B}^m$. Finally we look for
$\widetilde{u}$ which minimizes the energy defined as
$$
 \operatorname{argmin} _{u\in \mathbb{R}^{N_m},\,\Lambda\in
 \mathbb{R}^{N_\Lambda}} \frac{1}{2}\| u\| _{S_p}^2+\Lambda^\top \bigl(
 Cu-b\bigr), 
$$
As above, we take $\|u\|^2_{S_p} =\|(1-\Delta_\pi)^p u\|^2$. The problem
reduces to the regularized normal equation
\begin{equation}
  u=S_p^{-1}C^\top \bigl( CS_p^{-1}C^\top\bigr) ^{-1}b.
\end{equation} 
The linear system can then be solved using either the QR method or the
PCG method, as described in Section \ref{sec:method}. We remark that
the different extensions resulting from the different 
choice of $p$ used in the smoother will produce functions of a very
different nature, as they are minimizing different powers of the
Laplacian; which power $p$ is optimal will depend on the application
which the extension is being used for.

\begin{figure}
\includegraphics[scale=.4]{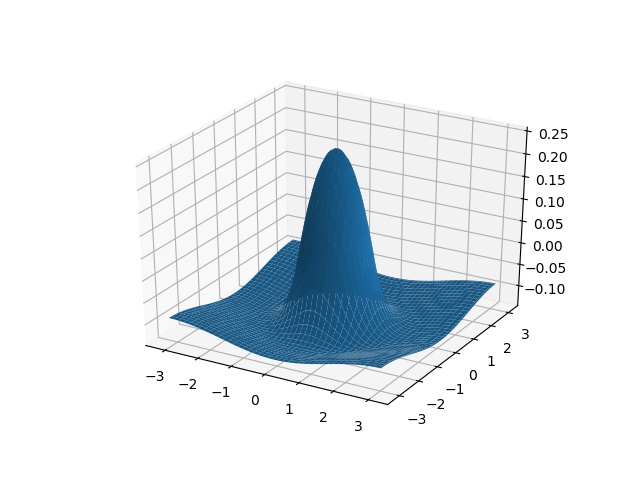}
\includegraphics[scale=.4]{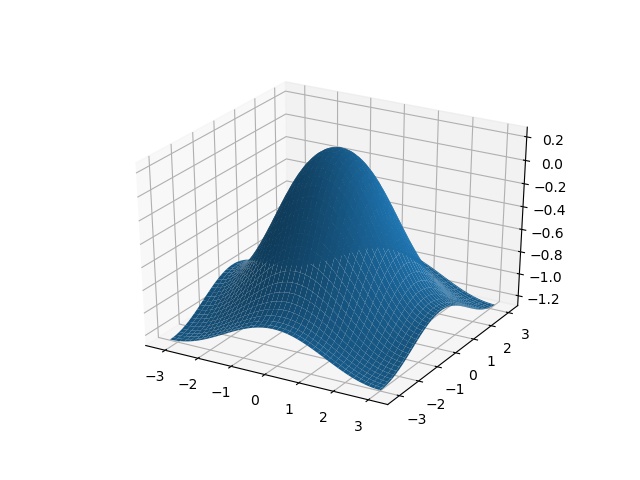}
\caption{$S_2$ and $S_4$ extensions of $\frac{1}{4}(1-r^2)$.}
\label{fig:extend}
\end{figure}

\begin{figure}
\includegraphics[scale=.4]{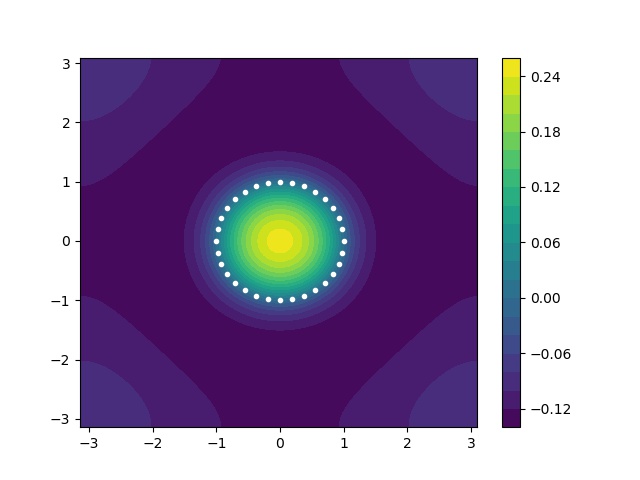}
\includegraphics[scale=.4]{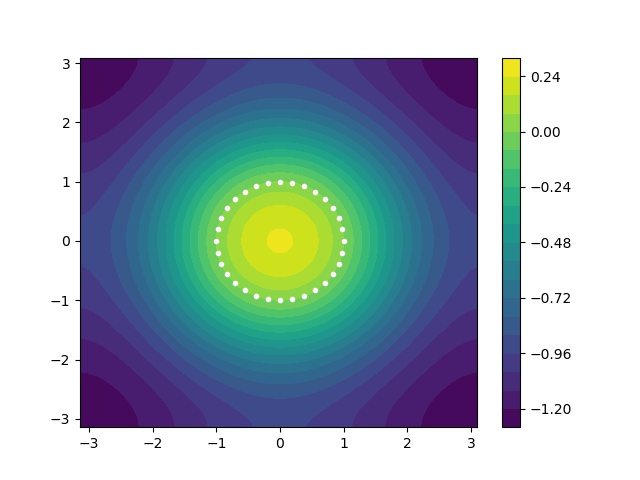}
\caption{$S_2$ and $S_4$ extensions of $\frac{1}{4}(1-r^2)$. The white
  dots show the discrete boundary between $\Omega$ and the fictitious
  domain.}
\label{fig:extendc}
\end{figure}

Next let's look at an example. Let $\Omega$ be the unit disc. We will
produce an extension of the function
$$u\big|_{\Omega} = \frac{1}{4}(1-r^2).$$
We do this by setting $S_pu =(1-\Delta_\pi)^p u$, where $p = 2,4$. The 
construction of $S_p^{-1}$ is the same as in Section
\ref{sec:regularizer}. The result of the extension done on a $128^2$ grid
using the two smoothers are shown in Figure \ref{fig:extend} with a
graph and in Figure \ref{fig:extendc} as a contour plot. In Table
\ref{tab:extend}, we show how the different smoothers affect
different Sobolev seminorms of the corresponding minimizers.
\begin{table}[H]
\centering
\begin{tabular}{|c|c|c|c|}
\hline
Smoother & $\|\nabla^2 u\|_{L_2}$ & $\|\nabla^3 u\|_{L_2}$ & $\|\nabla^4 u\|_{L_2}$ \\
\hline
$S_2$ & $1.99$ & $11.47$ & $169.24$\\
\hline
$S_4$ & $2.84$ & $3.42$ & $6.45$ \\
\hline
\end{tabular}
\caption{Gradient seminorms for the extension operator with $u =
  \frac{1}{4}(1-r^2)$ on a grid of size $128^2$.}
\label{tab:extend}
\end{table}

\section{Numerical Experiments}\label{sec:experiments}
In the numerical experiments of this section, we solve the
relevant system using both the QR and PCG methods following the
procedure described in Section \ref{sec:method}. We will then 
record the $L_2$ errors and the CPU times for each. We consider
problems with nonconstant coefficients, with complex
geometry, with Neumann boundary and with nonsmooth boundary
conditions. For last we solve a three dimensional problem with $\Omega
= B^3$, the ball of radius 1.
\subsection{Nonconstant Coefficients}\label{nonconstant}
\begin{figure}[H]
\includegraphics[scale=.45]{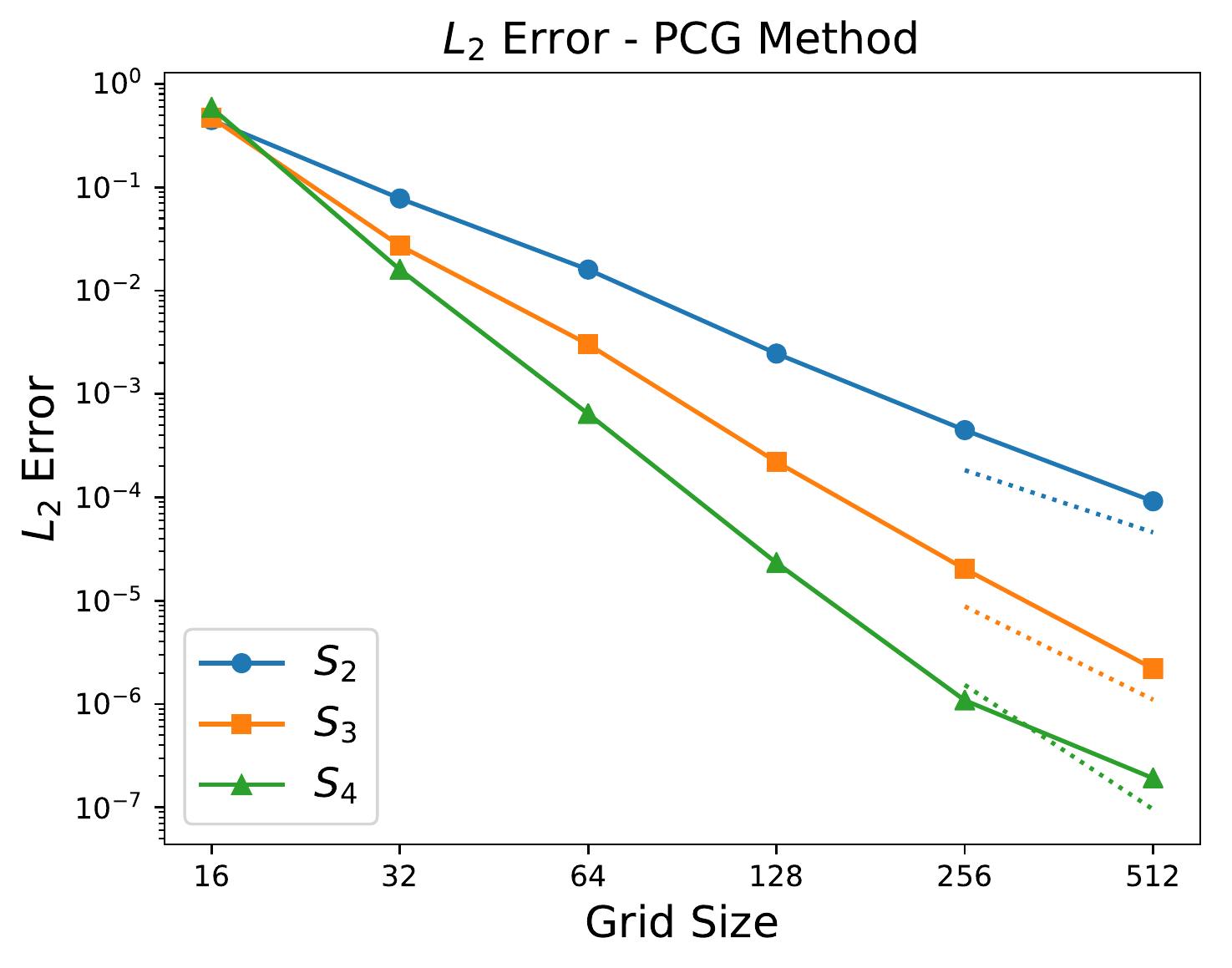}
\includegraphics[scale=.45]{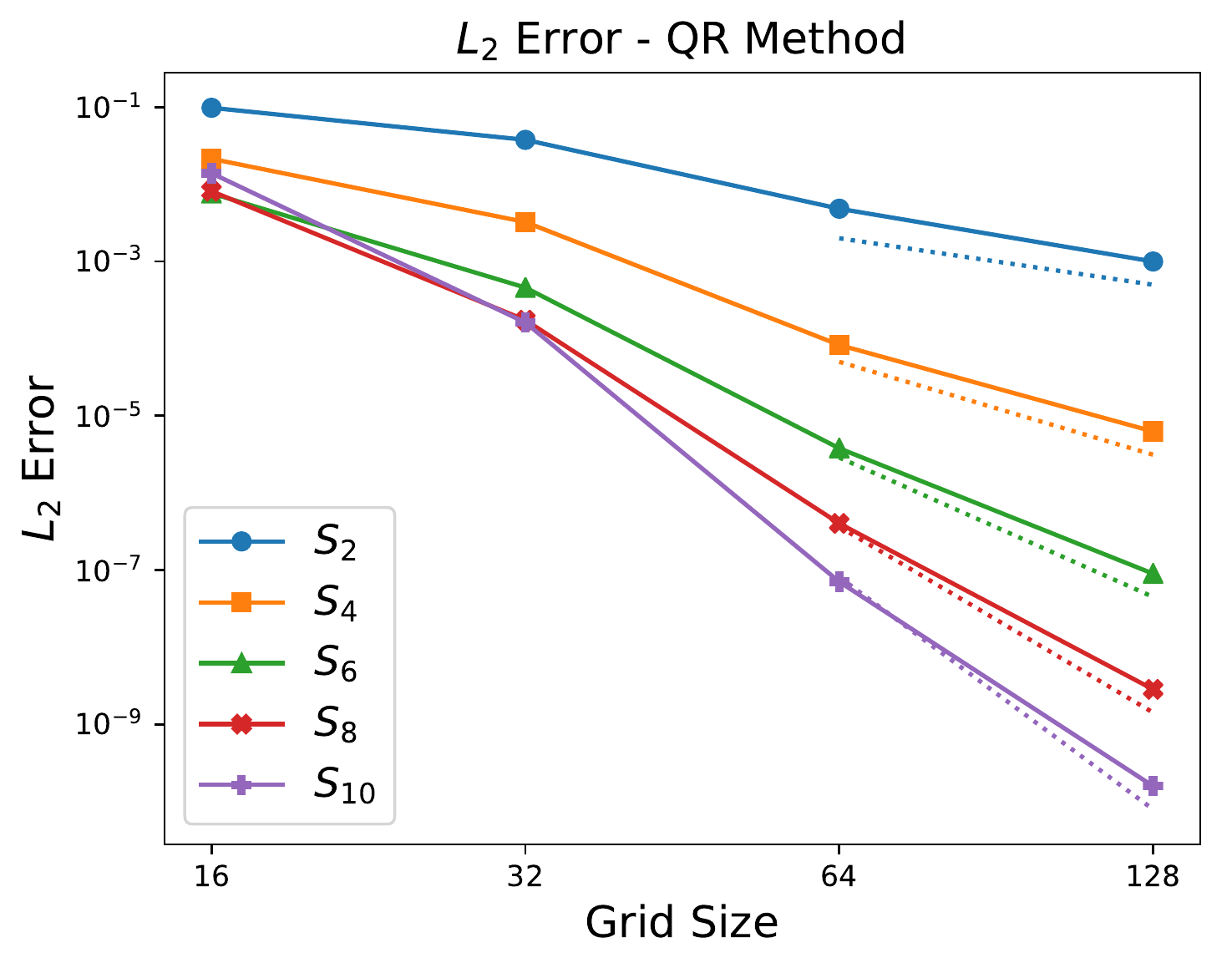}
\caption{Convergence of the $L_2$ error for different order smoothers solving Equation \protect\ref{eq:exp2}. The light dotted lines are reference lines of slope $\frac{1}{m^p}$ where $m$ is the number of grid points along one dimension.}
\label{fig:rates2}
\end{figure}

Let $\Omega$ be the unit disc $D$ discretized as
described in \ref{sec:domain} and study the problem 
\begin{equation}\label{eq:exp2}
\begin{cases}
-\Big[(2+y)\partial_x^2 + (2-x)\partial_y^2\Big] u = -6x(2+y) + 6y(2-x)&D, \\
u = x^3 - y^3&\text{on }\partial D.
\end{cases}
\end{equation}
The exact solution is $x^3 - y^3$.  The solution is calculated using the methodology described in Section \ref{sec:method}. Comparing Figures \ref{fig:rates1} and
\ref{fig:rates2}, we see that the accuracy achieved 
is roughly equivalent for both the constant coefficient problem and
the nonconstant coefficient problem. Comparing Table
\ref{tab:times2} with Table \ref{tab:times1}, we see that the QR method CPU time is comparable to the constant coefficient case, whereas the PCG method is considerably slower. Clearly, the condition number of the matrix grows faster for the nonconstant coefficient problem. However, both methods are still robust enough to efficiently solve nonconstant coefficient problems.

\begin{table}[H]
\centering
\begin{tabular}{|c|c|c|c|lllll}
\hline
\multicolumn{1}{|c|}{\multirow{2}{*}{\begin{tabular}[c]{@{}c@{}}Grid  \\ Size\end{tabular}}} & \multicolumn{3}{c|}{CPU Times - PCG Method} & \multicolumn{5}{c|}{CPU Times - QR Method} \\ \cline{2-9} 
\multicolumn{1}{|c|}{} & \multicolumn{1}{c|}{$S_2$} & \multicolumn{1}{c|}{$S_3$} & \multicolumn{1}{c|}{$S_4$} & \multicolumn{1}{c|}{$S_2$} & \multicolumn{1}{c|}{$S_4$} & \multicolumn{1}{c|}{$S_6$} & \multicolumn{1}{c|}{$S_8$} & \multicolumn{1}{c|}{$S_{10}$} \\ \hline
$16^2 $ & $0.24$ & $0.29$ & $0.29$ & \multicolumn{1}{c|}{$0.01$} & \multicolumn{1}{c|}{$0.01$} & \multicolumn{1}{c|}{$0.01$} & \multicolumn{1}{c|}{$0.01$} & \multicolumn{1}{c|}{$0.01$} \\
\hline
$32^2 $ & $0.27$ & $0.33$ & $0.39$ & \multicolumn{1}{c|}{$0.05$} & \multicolumn{1}{c|}{$0.05$} & \multicolumn{1}{c|}{$0.06$} & \multicolumn{1}{c|}{$0.06$} & \multicolumn{1}{c|}{$0.06$} \\
\hline
$64^2 $ & $0.3$ & $0.4$ & $0.46$ & \multicolumn{1}{c|}{$0.32$} & \multicolumn{1}{c|}{$0.45$} & \multicolumn{1}{c|}{$0.44$} & \multicolumn{1}{c|}{$0.47$} & \multicolumn{1}{c|}{$0.45$} \\
\hline
$128^2 $ & $0.52$ & $0.68$ & $0.93$ & \multicolumn{1}{c|}{$7.5$} & \multicolumn{1}{c|}{$9.0$} & \multicolumn{1}{c|}{$9.26$} & \multicolumn{1}{c|}{$9.41$} & \multicolumn{1}{c|}{$9.76$} \\
\hline
$256^2 $ & $1.95$ & $2.67$ & $5.94$& \multicolumn{5}{l}{\multirow{2}{*}{}} \\ \cline{1-4} 
$512^2 $ & $11.13$ & $14.51$ & $48.15$& \multicolumn{5}{l}{} \\ \cline{1-4}
\end{tabular}
\caption{CPU times for solving Equation \protect\ref{eq:exp2}. All computations were performed on an Intel 7700HQ.}
\label{tab:times2}
\end{table}

\subsection{A Flower Shaped Domain}\label{star}
\begin{figure}[H]
\includegraphics[scale=.45]{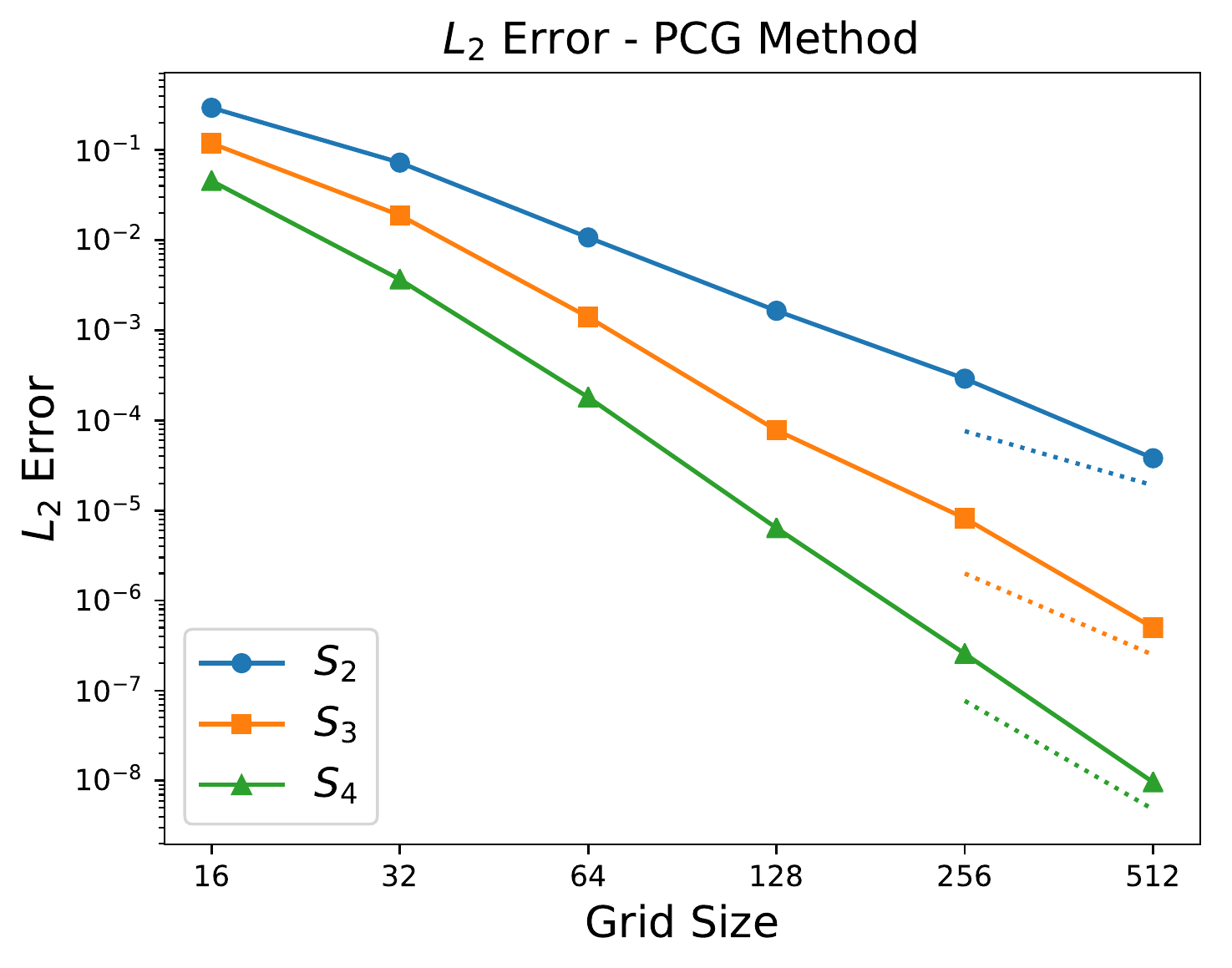}
\includegraphics[scale=.45]{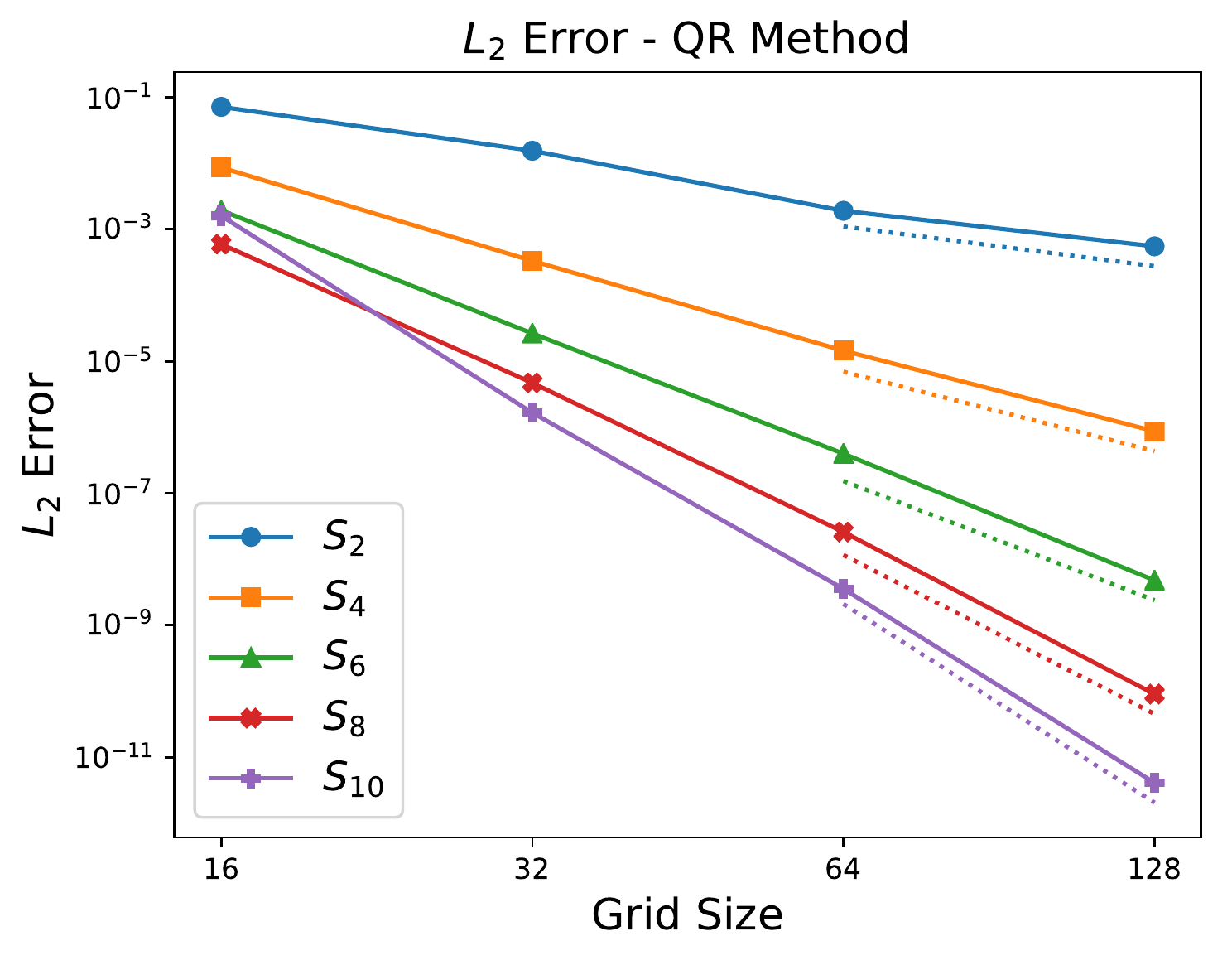}
\caption{Convergence of the $L_2$ error for different order smoothers solving Equation \protect\ref{eq:exp3}. The light dotted lines are reference lines of slope $\frac{1}{m^p}$ where $m$ is the number of grid points along one dimension.}
\label{fig:rates3}
\end{figure}

For a problem on a more complex domain, we consider a five petaled
flower. $$\Omega = \big\{(r,\theta) \;| \; r < 1 + .2\cos(5\theta)\big\}.$$
Figure \ref{fig:domains} shows $\Omega$ with its boundary discretization. We
solve the problem
\begin{equation}\label{eq:exp3}
\begin{cases}
-\Delta u = 0&\text{ in }\Omega,\\
\phantom{-\Delta}u= x^2 - y^2&\text{ on }
\partial \Omega
\end{cases}
\end{equation}
The exact solution is $x^2 - y^2$. The results of this experiment, contained in Figure \ref{fig:rates3} and Table \ref{tab:times3}, demonstrate that the accuracy and efficiency of the method is maintained even for a complex geometry. We note, however, that the PCG method is somewhat slower to converge on the flower than the disc.  We believe that, given the complexity of the shape, boundary points are necessarily closer together and therefore more
difficult for the interior grid to resolve. 

\begin{table}[H]
\centering
\begin{tabular}{|c|c|c|c|lllll}
\hline
\multicolumn{1}{|c|}{\multirow{2}{*}{\begin{tabular}[c]{@{}c@{}}Grid  \\ Size\end{tabular}}} & \multicolumn{3}{c|}{CPU Times - PCG Method} & \multicolumn{5}{c|}{CPU Times - QR Method} \\ \cline{2-9} 
\multicolumn{1}{|c|}{} & \multicolumn{1}{c|}{$S_2$} & \multicolumn{1}{c|}{$S_3$} & \multicolumn{1}{c|}{$S_4$} & \multicolumn{1}{c|}{$S_2$} & \multicolumn{1}{c|}{$S_4$} & \multicolumn{1}{c|}{$S_6$} & \multicolumn{1}{c|}{$S_8$} & \multicolumn{1}{c|}{$S_{10}$} \\ \hline
$16^2 $ & $0.27$ & $0.25$ & $0.26$ & \multicolumn{1}{c|}{$0.23$} & \multicolumn{1}{c|}{$0.21$} & \multicolumn{1}{c|}{$0.16$} & \multicolumn{1}{c|}{$0.17$} & \multicolumn{1}{c|}{$0.17$} \\
\hline
$32^2 $ & $0.29$ & $0.29$ & $0.32$ & \multicolumn{1}{c|}{$0.22$} & \multicolumn{1}{c|}{$0.2$} & \multicolumn{1}{c|}{$0.19$} & \multicolumn{1}{c|}{$0.22$} & \multicolumn{1}{c|}{$0.21$} \\
\hline
$64^2 $ & $0.37$ & $0.39$ & $0.46$ & \multicolumn{1}{c|}{$0.47$} & \multicolumn{1}{c|}{$0.54$} & \multicolumn{1}{c|}{$0.54$} & \multicolumn{1}{c|}{$0.54$} & \multicolumn{1}{c|}{$0.54$} \\
\hline
$128^2 $ & $0.51$ & $0.98$ & $1.56$ & \multicolumn{1}{c|}{$7.44$} & \multicolumn{1}{c|}{$8.34$} & \multicolumn{1}{c|}{$8.3$} & \multicolumn{1}{c|}{$8.29$} & \multicolumn{1}{c|}{$9.1$} \\
\hline
$256^2 $ & $1.51$ & $3.54$ & $16.48$& \multicolumn{5}{l}{\multirow{2}{*}{}} \\ \cline{1-4} 
$512^2 $ & $8.39$ & $23.78$ & $205.56$& \multicolumn{5}{l}{} \\ \cline{1-4}
\end{tabular}
\caption{CPU times for solving Equation \protect\ref{eq:exp3}. All computations were performed on an Intel 7700HQ.}
\label{tab:times3}
\end{table}

\subsection{Neumann Boundary Conditions}\label{sec:neumann}
\begin{figure}[H]
\includegraphics[scale=.45]{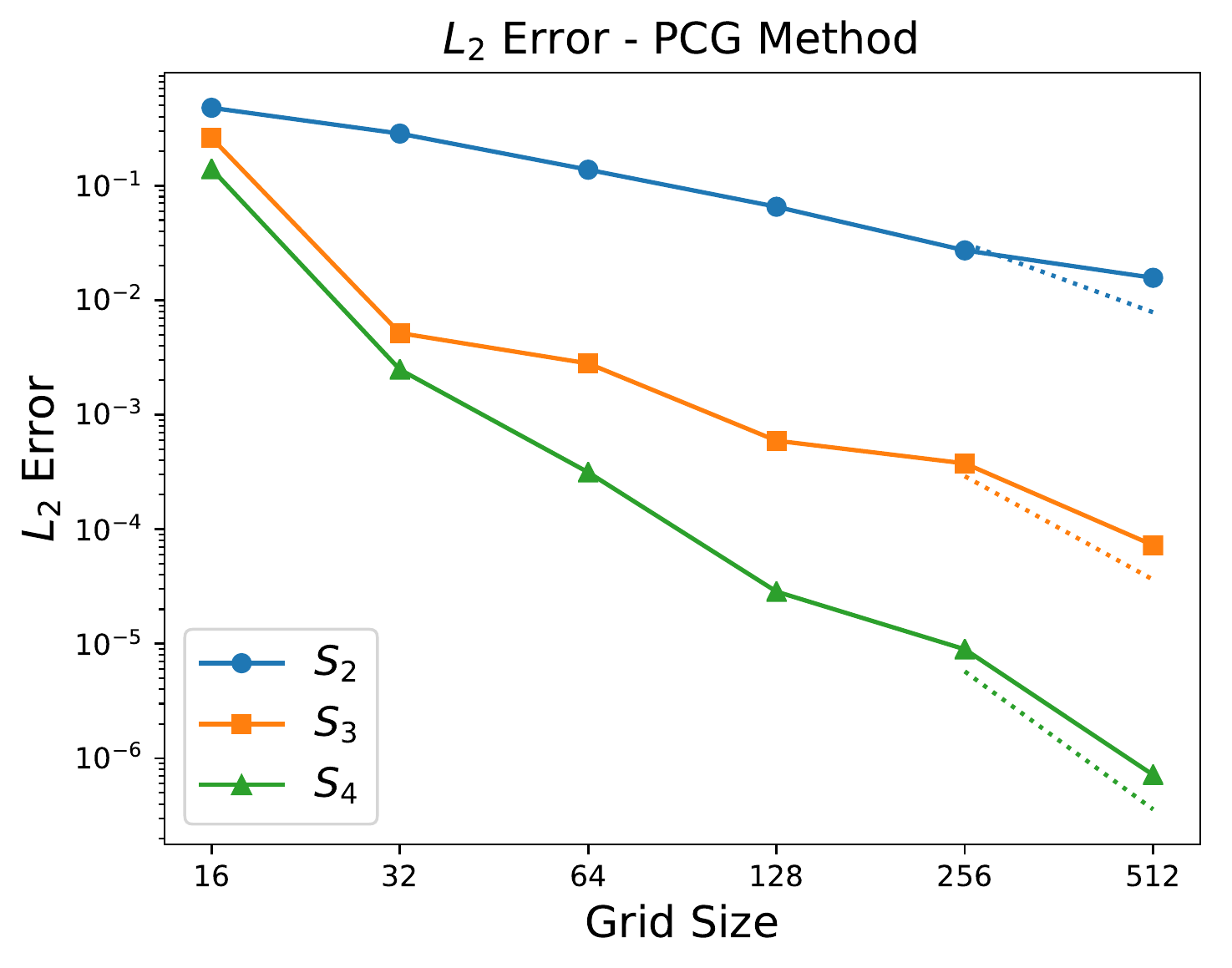}
\includegraphics[scale=.45]{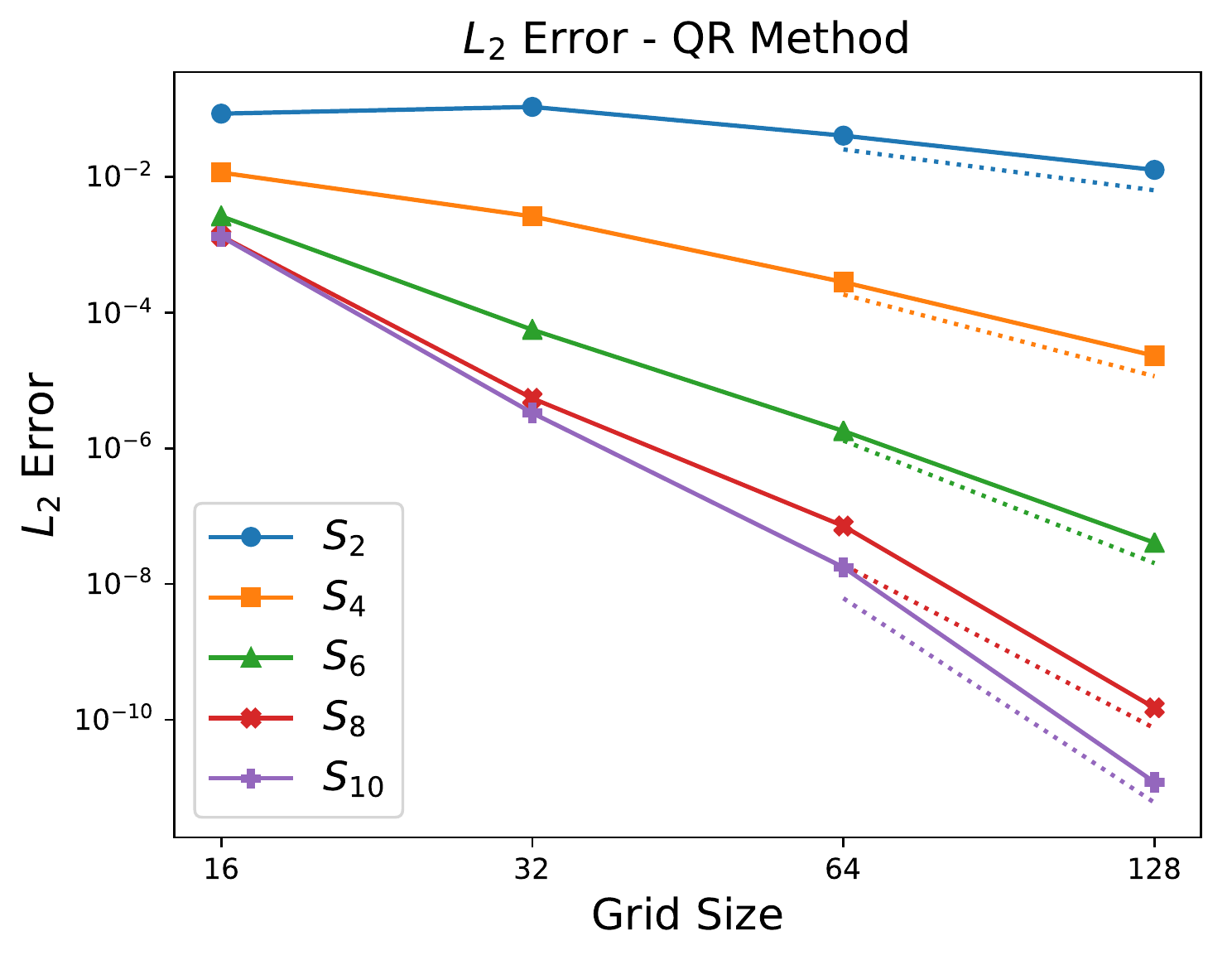}
\caption{Convergence of the $L_2$ error for different order smoothers solving Equation \protect\ref{eq:exp4}. The light dotted lines are reference lines of slope $\frac{1}{m^p}$ where $m$ is the number of grid points along one dimension.}
\label{fig:rates4}
\end{figure}

In order to demonstrate how the method is also applicable to other
boundary conditions, a Neumann problem is considered. We again set $\Omega$ as the unit disc. We solve the Neumann problem 
\begin{equation}\label{eq:exp4}
\begin{cases}
-\Delta u = 0 &\text{ in }\Omega,\\
\phantom{-\:}\frac{\partial u}{\partial \nu} = 
 2(x^2 - y^2)&\text{ on }\partial \Omega.
\end{cases}
\end{equation}
The exact solution is $x^2 - y^2$. As discussed in Section \ref{sec:method}, we evaluate the normal
derivative with a spectral interpolation. We substact $u(0,0)$ to eliminate the constant
functions in the kernel. The convergence results are displayed in
Figure \ref{fig:rates4}, and the CPU times are recorded in Table
\ref{tab:times4}. 

\begin{table}[H]
\centering
\begin{tabular}{|c|c|c|c|lllll}
\hline
\multicolumn{1}{|c|}{\multirow{2}{*}{\begin{tabular}[c]{@{}c@{}}Grid  \\ Size\end{tabular}}} & \multicolumn{3}{c|}{CPU Times - PCG Method} & \multicolumn{5}{c|}{CPU Times - QR Method} \\ \cline{2-9} 
\multicolumn{1}{|c|}{} & \multicolumn{1}{c|}{$S_2$} & \multicolumn{1}{c|}{$S_3$} & \multicolumn{1}{c|}{$S_4$} & \multicolumn{1}{c|}{$S_2$} & \multicolumn{1}{c|}{$S_4$} & \multicolumn{1}{c|}{$S_6$} & \multicolumn{1}{c|}{$S_8$} & \multicolumn{1}{c|}{$S_{10}$} \\ \hline
$16^2 $ & $0.22$ & $0.23$ & $0.22$ & \multicolumn{1}{c|}{$0.14$} & \multicolumn{1}{c|}{$0.16$} & \multicolumn{1}{c|}{$0.16$} & \multicolumn{1}{c|}{$0.17$} & \multicolumn{1}{c|}{$0.17$} \\
\hline
$32^2 $ & $0.23$ & $0.27$ & $0.27$ & \multicolumn{1}{c|}{$0.17$} & \multicolumn{1}{c|}{$0.23$} & \multicolumn{1}{c|}{$0.21$} & \multicolumn{1}{c|}{$0.2$} & \multicolumn{1}{c|}{$0.25$} \\
\hline
$64^2 $ & $0.26$ & $0.28$ & $0.31$ & \multicolumn{1}{c|}{$0.36$} & \multicolumn{1}{c|}{$0.56$} & \multicolumn{1}{c|}{$0.51$} & \multicolumn{1}{c|}{$0.5$} & \multicolumn{1}{c|}{$0.55$} \\
\hline
$128^2 $ & $0.42$ & $0.42$ & $0.61$ & \multicolumn{1}{c|}{$6.51$} & \multicolumn{1}{c|}{$8.15$} & \multicolumn{1}{c|}{$7.96$} & \multicolumn{1}{c|}{$8.14$} & \multicolumn{1}{c|}{$8.09$} \\
\hline
$256^2 $ & $1.23$ & $1.5$ & $3.29$& \multicolumn{5}{l}{\multirow{2}{*}{}} \\ \cline{1-4} 
$512^2 $ & $6.43$ & $8.32$ & $31.67$& \multicolumn{5}{l}{} \\ \cline{1-4}
\end{tabular}
\caption{CPU times for solving Equation \protect\ref{eq:exp4}. All computations were performed on an Intel 7700HQ.}
\label{tab:times4}
\end{table}

\subsection{Non-regular Problem}\label{sec:nonregular}
\begin{figure}[H]
\includegraphics[scale=.5]{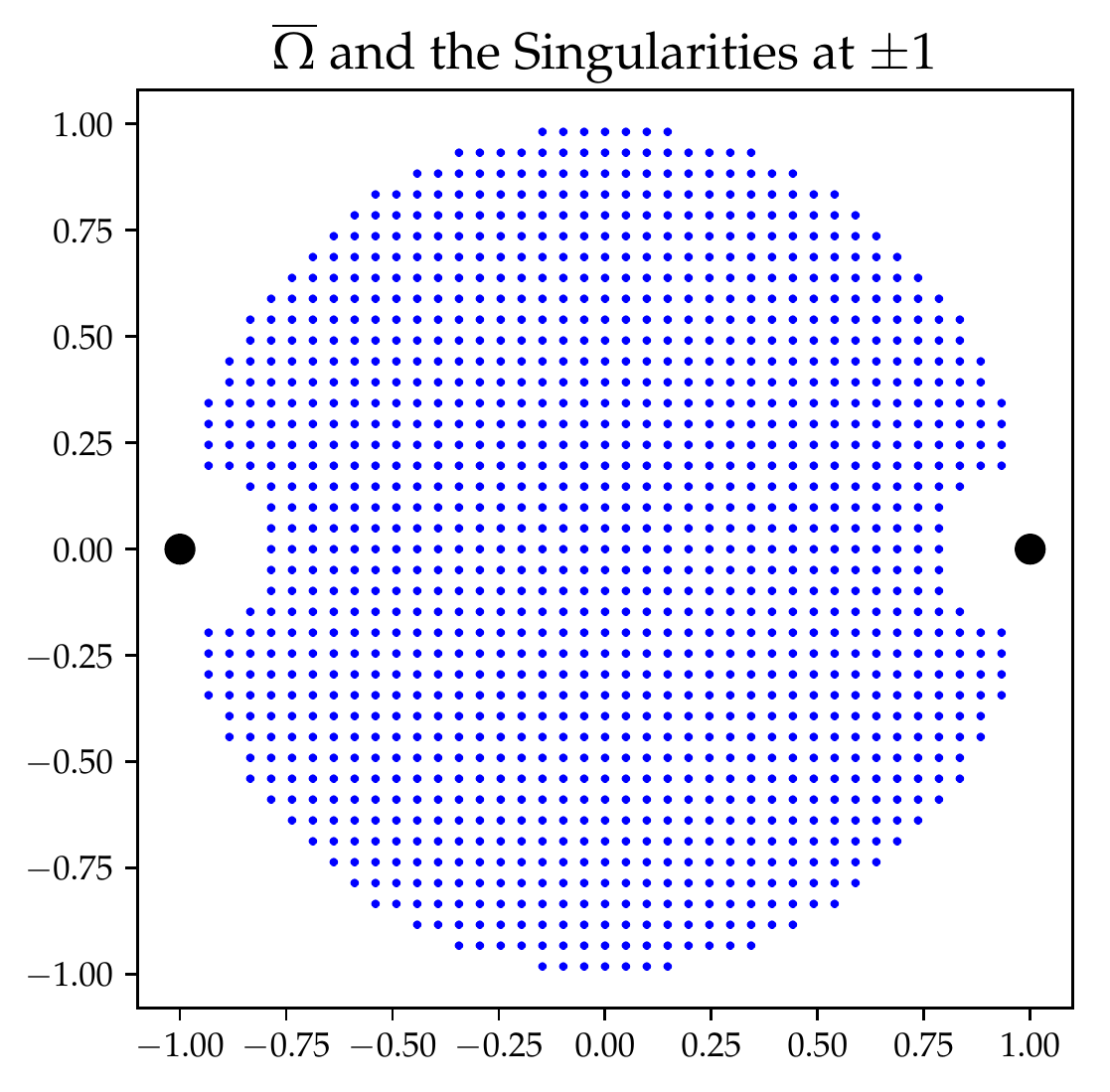}
\includegraphics[scale=.45]{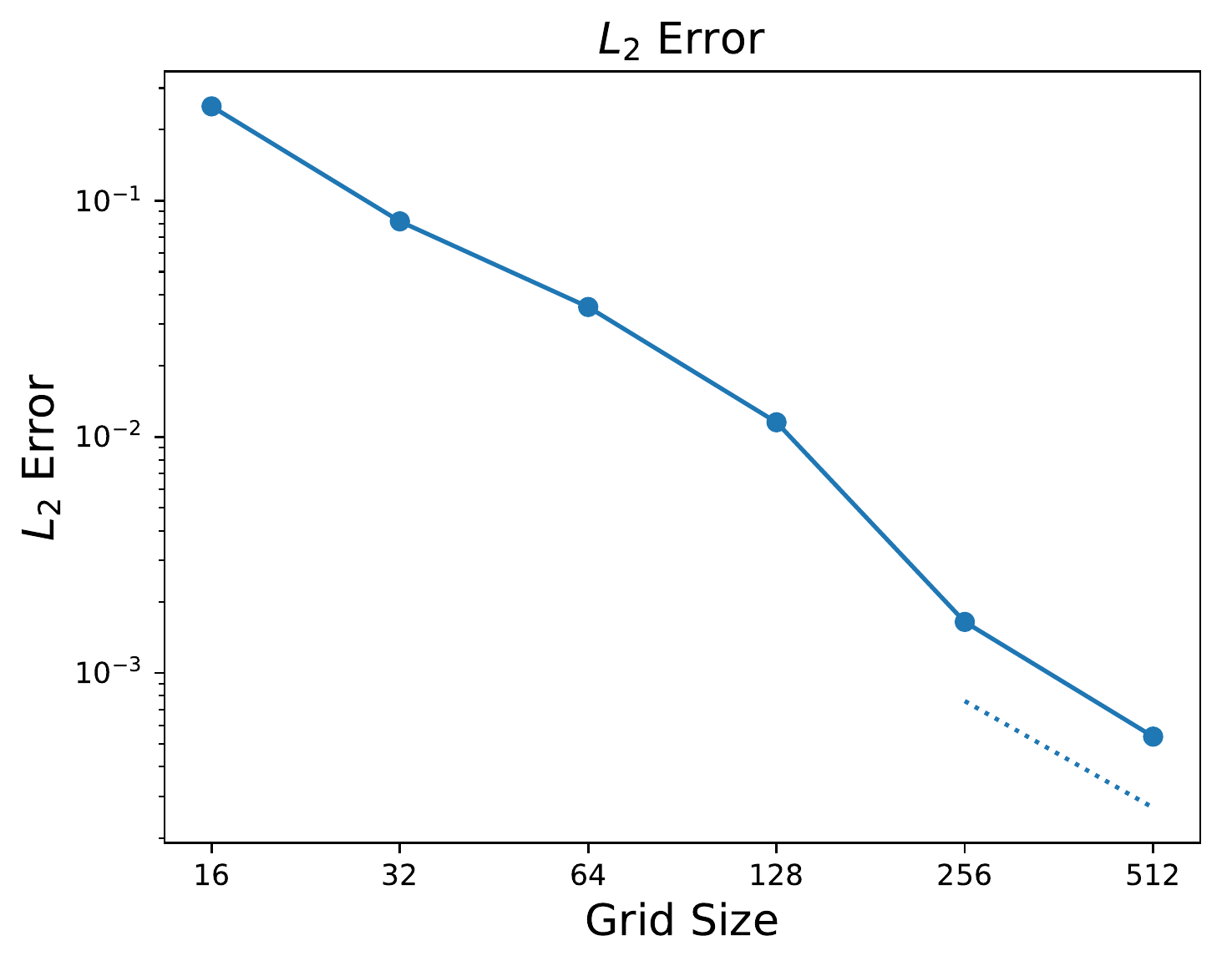}
\caption{The domain $\overline{\Omega}$ and the convergence of the
  nonregular problem away from the singularities. The reference line
  is slope $\left(\frac{1}{m}\right)^{1.5}$, where $m$ is the number
  of grid points along one dimension.} 
\label{fig:rates5}
\end{figure}

\begin{figure}[H]
\includegraphics[scale=.5]{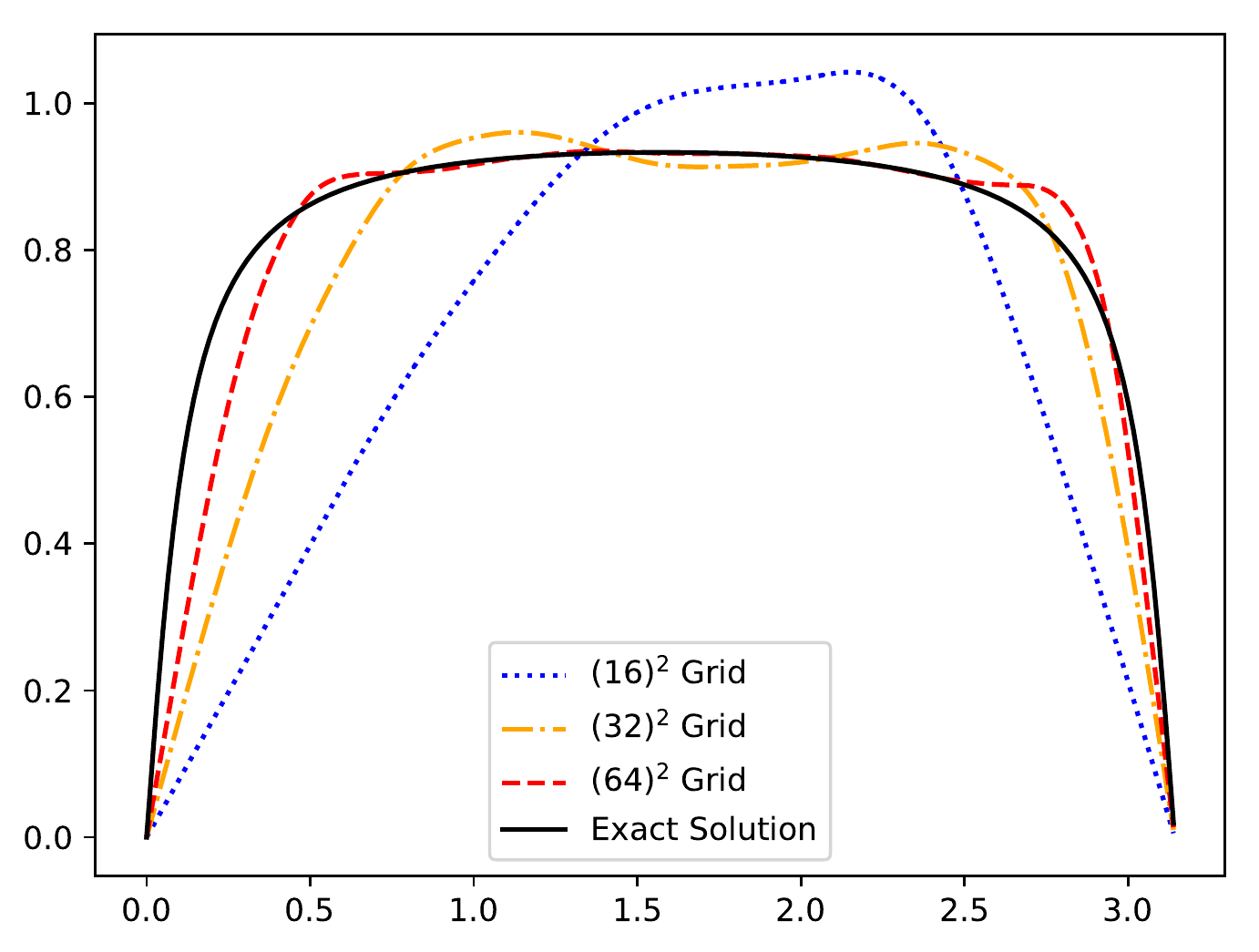}
\caption{Values of the approximated solution along the curve $r = .9$.}
\label{fig:level_sets5}
\end{figure}

Although our method is by its nature more suited to smooth problems,
it is not strictly limited to them. While the method in its current
form has no hope of properly approximating the solution in the
immediate vicinity of a singularity, outside a small ball containing
the singularity, it converges reasonably well to the solution. In the
following example, we let $\Omega$ be the unit disc and study the
nonsmooth problem 
$$\begin{cases}
-\Delta u = 0 &\text{ in }\Omega,\\
\phantom{-\Delta}u=g
&\text{ on }\partial \Omega,
\end{cases}$$
where
$$g(\theta)=\begin{cases}
1 &\text{ if } 0 \leq \theta \leq \pi \\
-1 &\text{ if } \pi \leq \theta \leq 2\pi
\end{cases}
$$
in polar coordinates. In solving this problem, we will use the same
discretization of the domain and operators used in the example studied
in Section \ref{sec:method}. We note that the higher order
smoothers do not provide an advantage when the solution itself is
not smooth, so we restrict ourselves to $S_2$. The true solution of this boundary value problem can be
given in the form of the series
$$u(r,\theta) = \sum_{k=1}^{\infty} g_k r^k \sin(k \theta),$$
where the Fourier coefficients $g_k$ are defined by
$$g_k = \begin{cases}
\frac{4}{k\pi} & \text{if k is odd,} \\
0 & \text{if k is even.}
\end{cases}$$
The two singularities occur at $y_1 = [1,0]$ and $y_2 = [-1,0]$. We
will study the solution away from the singularities in two
ways. First, we will look at the $L_2$ error on
$\widetilde{\Omega} = \Omega \backslash (B(y_1,0.2)\cup B(y_2,0.2))$,
which cuts out the singularities. The region $\widetilde{\Omega}$ and
the corresponding errors are shown in Figure
\ref{fig:rates5}. The graph show a convergence rate of approximately
$1.5$. We also show the approximated solution along 
the curve $r = .9$, $0 \leq \theta \leq \pi$ in Figure
\ref{fig:level_sets5}. As the grid becomes more dense, the
approximations get closer to the true solution.  

We also would like to point out that the general framework of our
method could potentially be modified to allow it to deal with singular
problems more effectively; this could either be done by allowing an
adaptive grid which is more dense in the region of the singularity or
by modifying the norm used to generate the smoother $S^{-1}$ by
introducing weights or allowing for some singular behavior. 

\subsection{The Sphere}\label{sec:sphere}
\begin{figure}[H]
\includegraphics[scale=.45]{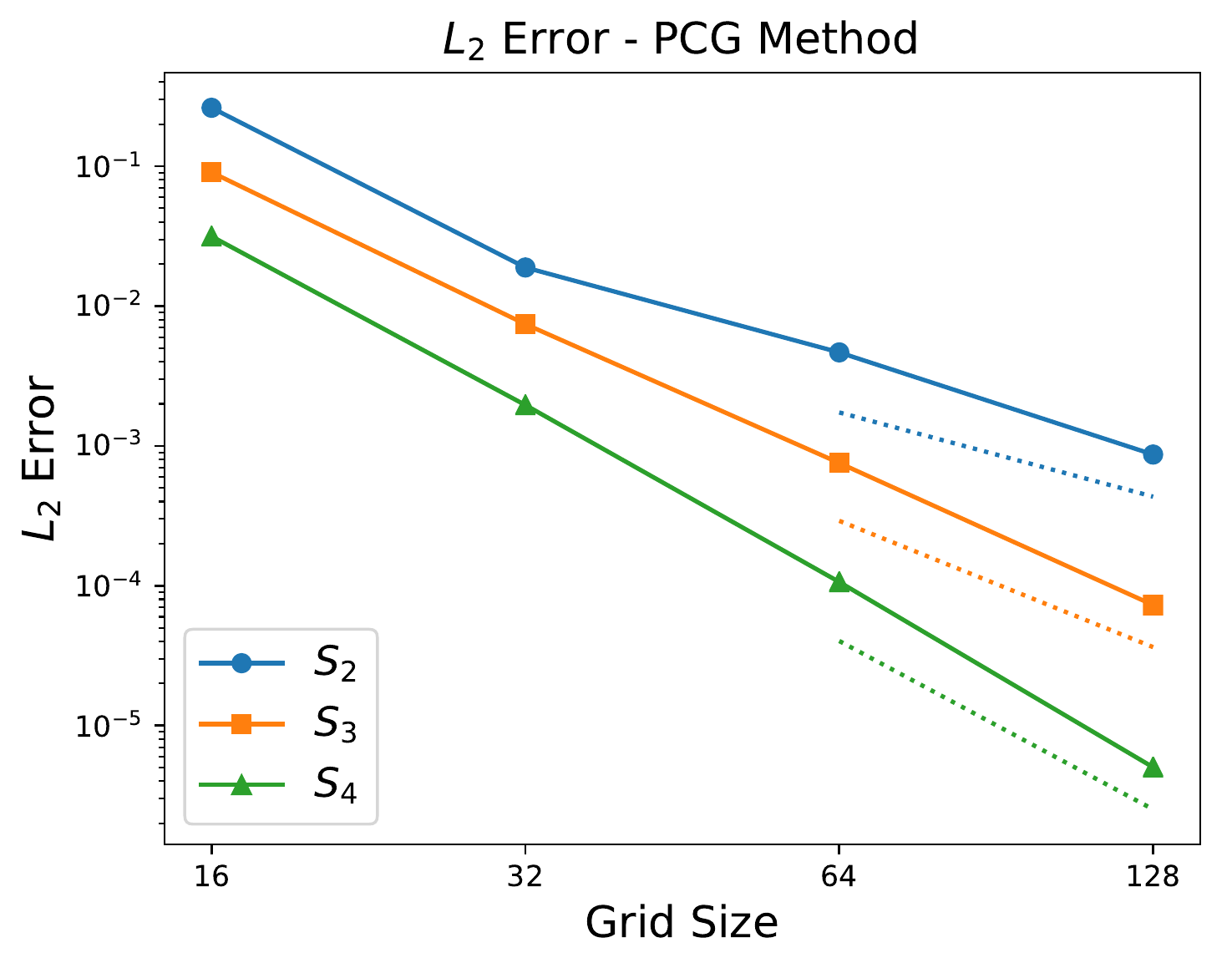}
\includegraphics[scale=.45]{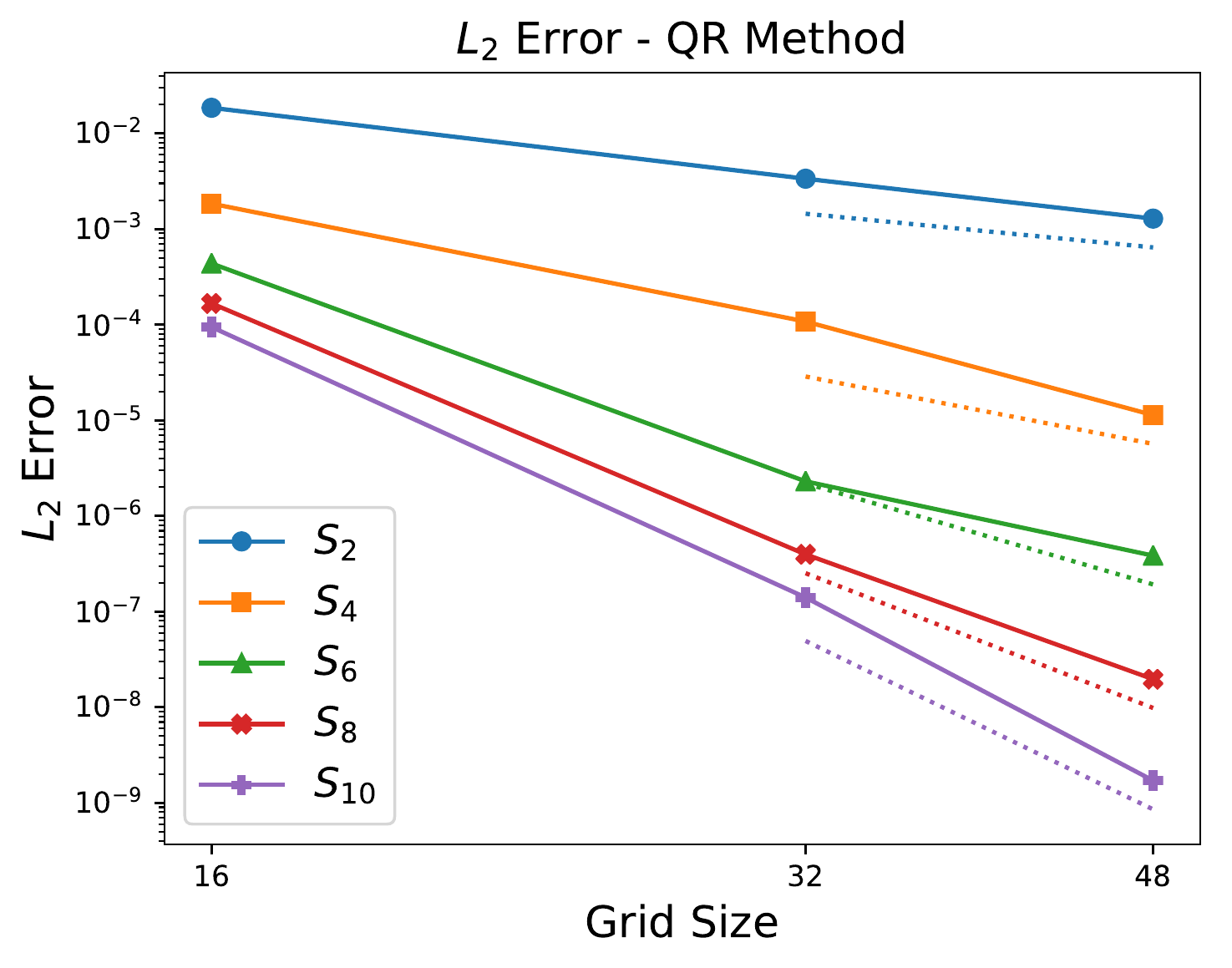}
\caption{Convergence of the $L_2$ error for different order smoothers
  solving Equation \protect\ref{eq:exp6}. The light dotted lines are
  reference lines of slope $\frac{1}{m^p}$ where $m$ is the number of
  grid points along one dimension.} 
\label{fig:rates6}
\end{figure}

For a three dimensional example, we choose the unit sphere embedded in
the three dimensional torus. The boundary is discretized by the well
known Fibonacci lattice \cite{FibSphere} which comes close to distributing
points uniformly on the sphere. With a box discretization of $m^3$
points, we use $\frac{1}{16}\left(\frac{m}{2\pi}\right)^2$ in the QR
method and $\frac{1}{16}\left(\frac{m}{2\pi}\right)^2$ boundary points
in the PCG method; see Section \ref{sec:domain} for a more detailed
discussion. The problem considered is
\begin{equation}\label{eq:exp6}
\begin{cases}
-\Delta u = 1 &\text{ in }\Omega,\\
\phantom{-\Delta}u=0&\text{ on }\partial \Omega.
\end{cases}
\end{equation}
The exact solution is $\frac{1-r^2}{6}$. We see in Figure
\ref{fig:rates6} that the rate of convergence achieved 
is similar to that of the two dimensional problem. Obviously, given
the larger dimension, the CPU times are significantly larger than in
two dimensions; however, the method is still quite fast. We note that
because we are calculating explicit matrices for the QR factorization,
RAM limitations prevented us from using a grid larger than $48^3$.

\begin{table}[H]
\centering
\begin{tabular}{|c|c|c|c|lllll}
\hline
\multicolumn{1}{|c|}{\multirow{2}{*}{\begin{tabular}[c]{@{}c@{}}Grid  \\ Size\end{tabular}}} & \multicolumn{3}{c|}{CPU Times - PCG Method} & \multicolumn{5}{c|}{CPU Times - QR Method} \\ \cline{2-9} 
\multicolumn{1}{|c|}{} & \multicolumn{1}{c|}{$S_2$} & \multicolumn{1}{c|}{$S_3$} & \multicolumn{1}{c|}{$S_4$} & \multicolumn{1}{c|}{$S_2$} & \multicolumn{1}{c|}{$S_4$} & \multicolumn{1}{c|}{$S_6$} & \multicolumn{1}{c|}{$S_8$} & \multicolumn{1}{c|}{$S_{10}$} \\ \hline
$16^3 $ & $0.02$ & $0.12$ & $0.06$ & \multicolumn{1}{c|}{$0.23$} & \multicolumn{1}{c|}{$0.25$} & \multicolumn{1}{c|}{$0.35$} & \multicolumn{1}{c|}{$0.29$} & \multicolumn{1}{c|}{$0.29$} \\
\hline
$32^3 $ & $0.09$ & $0.16$ & $0.33$ & \multicolumn{1}{c|}{$9.08$} & \multicolumn{1}{c|}{$12.07$} & \multicolumn{1}{c|}{$11.04$} & \multicolumn{1}{c|}{$11.51$} & \multicolumn{1}{c|}{$11.5$} \\
\hline
$48^3 $ & $0.44$ & $0.70$ & $1.44$ & \multicolumn{1}{c|}{$189.78$} & \multicolumn{1}{c|}{$196.58$} & \multicolumn{1}{c|}{$175.58$} & \multicolumn{1}{c|}{$182.03$} & \multicolumn{1}{c|}{$185.81$} \\
\hline
$64^3 $ & $1.51$ & $2.81$ & $5.96$& \multicolumn{5}{l}{\multirow{2}{*}{}} \\ \cline{1-4} 
$128^3 $ & $11.54$ & $22.11$ & $161.78$& \multicolumn{5}{l}{} \\ \cline{1-4}
\end{tabular}
\caption{CPU times for solving Equation \protect\ref{eq:exp6}. All computations were performed on an Intel 7700HQ.}
\label{tab:times6}
\end{table}

\end{document}